\input amssym.tex
\input cyracc.def

\magnification =1000
\vsize = 195 true mm
\hsize =125 true mm
\baselineskip=12truept
\voffset=2\baselineskip
\overfullrule = 0pt
\baselineskip=15truept
\parindent=0pt

\pageno=1

\headline= {\ifodd\pageno\rightheadline \else\leftheadline\fi}
\def\rightheadline{\tenrm\hfil Equivalent Inequalities \hfil\folio}
\def\leftheadline{\tenrm\folio\hfil Equivalent Inequalities\hfil}


\def\func#1#2#3{$#1\colon #2\mapsto #3$}

\def\equivi{\! \equiv\! }
\def\section{\sectn}
\def\wide#1{\vskip 8pt\baselineskip=14truept \parshape= 8 40pt10cm40pt10cm40pt10cm40pt10cm40pt10cm40pt10cm40pt10cm40pt10cm\ {\sl Convention:\  #1}\bigskip\baselineskip
15truept}
\def\fa{\forall\, }
\def\we{\wedge\,}
\def\gwide#1{\vskip 4pt\baselineskip=14truept \parshape= 8 40pt10cm40pt10cm40pt10cm40pt10cm40pt10cm40pt10cm40pt10cm40pt10cm\ {\sl  \ #1}\vskip4pt\baselineskip
15truept}

\def \BR{\Bbb R}
\def \BN{\Bbb N}
\def\BP{\Bbb P}

\def \BZ{\Bbb Z}

\def\Lra{\Longrightarrow}
\def\and{\wedge\,}
\def\or{\vee\,}

\def\bvec#1{{\bf #1}}

 \def\recip#1{\frac{1}{#1}}

 \def\GA{{\rm (\!GA\!)}}

 \def\frac#1#2{{#1\over #2}}
 \def\dfrac{\displaystyle\frac}
 
 \def\noter#1{\advance\notenumber by 1
$^{\rn\the\notenumber}$\footnote{}{$^{\rn\the\notenumber}$\sevenrm #1}}
 \def\inv#1{#1^{-1}}
 \def\refb#1{[{\it #1\/}]}
 \def\rn{\romannumeral}
 
 \def\ref#1{[{\sl #1\/}]}
 \def\qed{{$\square$}\qquad}%
 \def\QED{\hfill{$\square$}\smallskip}%

 \outer\def\proclaim#1#2{\par\medbreak\noindent{\bf #1 \enspace}{\sl #2}\par\ifdim\lastskip<\medskipamount
 \removelastskip\penalty55\medskip\fi}%
 \outer\def\proclaimsc#1#2{\par\medbreak\noindent{\ninesc #1 \enspace}{\sl
 #2}\par\ifdim\lastskip<\medskipamount \removelastskip\penalty55\medskip\fi}%
 \outer\def\proclaimbf#1#2{\par\medbreak\noindent{\bfn #1 \enspace}{\sl
 #2}\par\ifdim\lastskip<\medskipamount \removelastskip\penalty55\medskip\fi}%
 \outer\def\proclaimbfn#1#2{\par\medbreak\noindent{\bfnn #1 \enspace}{\sl
 #2}\par\ifdim\lastskip<\medskipamount \removelastskip\penalty55\medskip\fi}%

%

%

\newfam\cyrfam
\font\fivecyr=wncyr5
\font\sevencyr=wncyr10 at 7truept

\font\tencyr=wncyr10

\def\cyrseven{\sevencyr\cyracc}

\def\cyrf{\fam\cyrfam\sevencyr\cyracc}

\textfont\cyrfam=\tencyr 
\scriptfont\cyrfam=\sevencyr
\scriptscriptfont\cyrfam=\fivecyr

\font\ssbx=cmssbx10 at 20pt 
\font\teneufm=eufm10

\def\eufm#1{{\hbox{\teneufm #1}}}

\font\sevenit=cmti7


\font\bfn=cmbx9
\font\bfnn=cmbx8
\font\bfnnn=cmbx8 at 7pt
\font\bfx=cmbx12

\font\ninesc=cmcsc9
\font\ssbx=cmssbx10 at 20pt 
\def\cal{\fam2\tensy}


\def\chap#1#2#3{\null\vskip 2truecm \noindent{{\ssbx
\hbox{#1\quad  #2}\smallskip\hbox{ \hskip 1cm  #3}} }\vskip
1.3truecm}

\def\abstract#1{\baselineskip=5pt \parshape= 8
40pt8cm40pt8cm40pt8cm40pt8cm40pt8cm40pt8cm 40pt8cm40pt8cm{\bfn  Abstract}:{\sevenrm #1}\bigskip\baselineskip
15pt}

\newcount\notenumber

%
%

\newcount\sectnumber%
\def\sectn#1 {\advance \sectnumber by 1{\par\bigskip\noindent{\bfx
\the\sectnumber\ #1}}\par}%

\newcount\subsanumber%
\def\asubsn#1 {\advance \subsanumber by 1{\par\medskip\noindent{\bf\the\sectnumber.\the\subsanumber\ 
#1}\quad}}%

\newcount\subsbnumber%
\def\bsubsn#1 {\advance \subsbnumber by 1{\par\medskip\noindent{\bf\the\sectnumber.\the\subsbnumber\ 
#1}\quad}}%

\newcount\subscnumber%
\def\csubsn#1 {\advance \subscnumber by 1{\par\medskip\noindent{\bf\the\sectnumber.\the\subscnumber\ 
#1}\quad}}%

\newcount\subsdnumber%
\def\dsubsn#1 {\advance \subsdnumber by 1{\par\medskip\noindent{\bf\the\sectnumber.\the\subsdnumber\ 
#1}\quad}}%

\newcount\subsenumber%
\def\esubsn#1 {\advance \subsenumber by 1{\par\medskip\noindent{\bf\the\sectnumber.\the\subsenumber\ 
#1}\quad}}

\newcount\subsfnumber%
\def\fsubsn#1 {\advance \subsfnumber by 1{\par\medskip\noindent{\bf\the\sectnumber.\the\subsfnumber\ 
#1}\quad}}%


\newcount\subssaanumber%
\def\aasubssn#1 {\advance \subssaanumber by
1{\par\medskip\noindent{\bfn \the\sectnumber.\the\subsanumber.\the\subssaanumber\  #1}\quad}}%

\newcount\subssabnumber%
\def\absubssn#1 {\advance \subssabnumber by
1{\par\medskip\noindent{\bfn \the\sectnumber.\the\subsanumber.\the\subssabnumber\  #1}\quad}}%

\newcount\subssbanumber%
\def\basubssn#1 {\advance \subssbanumber by
1{\par\medskip\noindent{\bfn \the\sectnumber.\the\subsbnumber.\the\subssbanumber\  #1}\quad}}%

\newcount\subssbbnumber%
\def\bbsubssn#1 {\advance \subssbbnumber by
1{\par\medskip\noindent{\bfn \the\sectnumber.\the\subsbnumber.\the\subssbbnumber\  #1}\quad}}%

\newcount\subssbcnumber%
\def\bcsubssn#1 {\advance \subssbcnumber by
1{\par\medskip\noindent{\bfn \the\sectnumber.\the\subsbnumber.\the\subssbcnumber\  #1}\quad}}%

\newcount\subssbdnumber%
\def\bdsubssn#1 {\advance \subssbdnumber by
1{\par\medskip\noindent{\bfn \the\sectnumber.\the\subsbnumber.\the\subssbdnumber\  #1}\quad}}%

\newcount\subssbenumber%
\def\besubssn#1 {\advance \subssbenumber by
1{\par\medskip\noindent{\bfn \the\sectnumber.\the\subsbnumber.\the\subssbenumber\  #1}\quad}}%

\newcount\subsscanumber%
\def\casubssn#1 {\advance \subsscanumber by
1{\par\medskip\noindent{\bfn \the\sectnumber.\the\subscnumber.\the\subsscanumber\  #1}\quad}}\relax%

\newcount\subsscbnumber%
\def\cbsubssn#1 {\advance \subsscbnumber by
1{\par\medskip\noindent{\bfn \the\sectnumber.\the\subscnumber.\the\subsscbnumber\  #1}\quad}}%

\newcount\subssccnumber%
\def\ccsubssn#1 {\advance \subssccnumber by
1{\par\medskip\noindent{\bfn \the\sectnumber.\the\subscnumber.\the\subssccnumber\  #1}\quad}}%

\newcount\subssdanumber%
\def\dasubssn#1 {\advance \subssdanumber by
1{\par\medskip\noindent{\bfn \the\sectnumber.\the\subsdnumber.\the\subssdanumber\  #1}\quad}}%

\newcount\subssdbnumber%
\def\dbsubssn#1 {\advance \subssdbnumber by
1{\par\medskip\noindent{\bfn \the\sectnumber.\the\subsdnumber.\the\subssdbnumber\  #1}\quad}}%

\newcount\subssdcnumber%
\def\dcsubssn#1 {\advance \subssdcnumber by
1{\par\medskip\noindent{\bfn \the\sectnumber.\the\subsdnumber.\the\subssdcnumber\  #1}\quad}}%

\newcount\subsseanumber%
\def\easubssn#1 {\advance \subsseanumber by
1{\par\medskip\noindent{\bfn \the\sectnumber.\the\subsenumber.\the\subsseanumber\  #1}\quad}}%

\newcount\subssebnumber%
\def\ebsubssn#1 {\advance \subssebnumber by
1{\par\medskip\noindent{\bfn \the\sectnumber.\the\subsenumber.\the\subssebnumber\  #1}\quad}}%

\newcount\subssecnumber%
\def\ecsubssn#1 {\advance \subssecnumber by
1{\par\medskip\noindent{\bfn \the\sectnumber.\the\subsenumber.\the\subssecnumber\  #1}\quad}}%

\newcount\subssefnumber%
\def\efsubssn#1 {\advance \subssefnumber by
1{\par\medskip\noindent{\bfn \the\sectnumber.\the\subsenumber.\the\subssefnumber\  #1}\quad}}%

\newcount\subssehnumber%
\def\ehsubssn#1 {\advance \subssehnumber by
1{\par\medskip\noindent{\bfn \the\sectnumber.\the\subsenumber.\the\subssehnumber\  #1}\quad}}%

\newcount\subssffnumber%
\def\ffsubssn#1 {\advance \subssffnumber by
1{\par\medskip\noindent{\bfn \the\sectnumber.\the\subsfnumber.\the\subssffnumber\  #1}\quad}}%

\newcount\subssfhnumber%
\def\fhsubssn#1 {\advance \subssfhnumber by
1{\par\medskip\noindent{\bfn \the\sectnumber.\the\subsfnumber.\the\subssfhnumber\  #1}\quad}}%

\newcount\subsssaaanumber%
\def\aaasubsssn#1 {\advance \subsssaaanumber by
1{\par\medskip\noindent{\bfnn
\the\sectnumber.\the\subsanumber.\the\subssaanumber.\the\subsssaaanumber\ #1}\quad}}%

\newcount\subsssabanumber%
\def\abasubsssn#1 {\advance \subsssabanumber by
1{\par\medskip\noindent{\bfnn
\the\sectnumber.\the\subsanumber.\the\subssabnumber.\the\subsssabanumber\ #1}\quad}}%

\newcount\subsssabbnumber%
\def\abbsubsssn#1 {\advance \subsssabbnumber by
1{\par\medskip\noindent{\bfnn
\the\sectnumber.\the\subsanumber.\the\subssabnumber.\the\subsssabbnumber\ #1}\quad}}%

\newcount\subsssabcnumber%
\def\abcsubsssn#1 {\advance \subsssabcnumber by
1{\par\medskip\noindent{\bfnn
\the\sectnumber.\the\subsanumber.\the\subssabnumber.\the\subsssabcnumber\ #1}\quad}}%

\newcount\subsssabdnumber%
\def\abdsubsssn#1 {\advance \subsssabdnumber by
1{\par\medskip\noindent{\bfnn
\the\sectnumber.\the\subsanumber.\the\subssabnumber.\the\subsssabdnumber\
#1}\quad}}%

\newcount\subsssabenumber%
\def\abesubsssn#1 {\advance \subsssabenumber by
1{\par\medskip\noindent{\bfnn
\the\sectnumber.\the\subsanumber.\the\subssabnumber.\the\subsssabenumber\
#1}\quad}}%

\newcount\subsssacanumber%
\def\acasubsssn#1 {\advance \subsssacanumber by
1{\par\medskip\noindent{\bfnn
\the\sectnumber.\the\subsanumber.\the\subssacnumber.\the\subsssacanumber\ #1}\quad}}%

\newcount\subsssbaanumber%
\def\baasubsssn#1 {\advance \subsssbaanumber by
1{\par\medskip\noindent{\bfnn \the\sectnumber.\the\subsbnumber.\the\subssbanumber.\the\subsssbaanumber\ #1}\quad}}%

\newcount\subsssbacnumber%
\def\bacsubsssn#1 {\advance \subsssbacnumber by
1{\par\medskip\noindent{\bfnn \the\sectnumber.\the\subsbnumber.\the\subssbanumber.\the\subsssbacnumber\ #1}\quad}}%

\newcount\subsssbadnumber%
\def\badsubsssn#1 {\advance \subsssbadnumber by
1{\par\medskip\noindent{\bfnn \the\sectnumber.\the\subsbnumber.\the\subssbanumber.\the\subsssbadnumber\ #1}\quad}}%

\newcount\subsssbaenumber%
\def\baesubsssn#1 {\advance \subsssbaenumber by
1{\par\medskip\noindent{\bfnn \the\sectnumber.\the\subsbnumber.\the\subssbanumber.\the\subsssbaenumber\ #1}\quad}}%

\newcount\subsssbbanumber%
\def\bbasubsssn#1 {\advance \subsssbbanumber by
1{\par\medskip\noindent{\bfnn \the\sectnumber.\the\subsbnumber.\the\subssbbnumber.\the\subsssbbanumber\ #1}\quad}}%

\newcount\subsssbbbnumber%
\def\bbbsubsssn#1 {\advance \subsssbbbnumber by
1{\par\medskip\noindent{\bfnn \the\sectnumber.\the\subsbnumber.\the\subssbbnumber.\the\subsssbbbnumber\ #1}\quad}}%

\newcount\subsssbbcnumber%
\def\bbcsubsssn#1 {\advance \subsssbbcnumber by
1{\par\medskip\noindent{\bfnn \the\sectnumber.\the\subsbnumber.\the\subssbbnumber.\the\subsssbbcnumber\ #1}\quad}}%

\newcount\subsssbbdnumber%
\def\bbdsubsssn#1 {\advance \subsssbbdnumber by
1{\par\medskip\noindent{\bfnn \the\sectnumber.\the\subsbnumber.\the\subssbbnumber.\the\subsssbbdnumber\ #1}\quad}}%

\newcount\subsssbcanumber%
\def\bcasubsssn#1 {\advance \subsssbcanumber by
1{\par\medskip\noindent{\bfnn \the\sectnumber.\the\subsbnumber.\the\subssbcnumber.\the\subsssbcanumber\ #1}\quad}}%

\newcount\subsssbccnumber%
\def\bccsubsssn#1 {\advance \subsssbccnumber by
1{\par\medskip\noindent{\bfnn \the\sectnumber.\the\subsbnumber.\the\subssbcnumber.\the\subsssbccnumber\ #1}\quad}}%

\newcount\subsssbdanumber%
\def\bdasubsssn#1 {\advance \subsssbdanumber by
1{\par\medskip\noindent{\bfnn \the\sectnumber.\the\subsbnumber.\the\subssbdnumber.\the\subsssbdanumber\ #1}\quad}}%

\newcount\subsssbdbnumber%
\def\bdbsubsssn#1 {\advance \subsssbdbnumber by
1{\par\medskip\noindent{\bfnn \the\sectnumber.\the\subsbnumber.\the\subssbdnumber.\the\subsssbdbnumber\ #1}\quad}}%

\newcount\subsssbeanumber%
\def\beasubsssn#1 {\advance \subsssbeanumber by
1{\par\medskip\noindent{\bfnn \the\sectnumber.\the\subsbnumber.\the\subssbenumber.\the\subsssbeanumber\ #1}\quad}}%

\newcount\subsssbebnumber%
\def\bebsubsssn#1 {\advance \subsssbebnumber by
1{\par\medskip\noindent{\bfnn \the\sectnumber.\the\subsbnumber.\the\subssbenumber.\the\subsssbebnumber\ #1}\quad}}%

\newcount\subsssbecnumber%
\def\becsubsssn#1 {\advance \subsssbecnumber by
1{\par\medskip\noindent{\bfnn \the\sectnumber.\the\subsbnumber.\the\subssbenumber.\the\subsssbecnumber\ #1}\quad}}%

\newcount\subsssbednumber%
\def\bedsubsssn#1 {\advance \subsssbednumber by
1{\par\medskip\noindent{\bfnn \the\sectnumber.\the\subsbnumber.\the\subssbenumber.\the\subsssbednumber\ #1}\quad}}%

\newcount\subssscaanumber%
\def\caasubsssn#1 {\advance \subssscaanumber by
1{\par\medskip\noindent{\bfnn \the\sectnumber.\the\subscnumber.\the\subsscanumber.\the\subssscaanumber\ #1}\quad}}%

\newcount\subssscacnumber%
\def\cacsubsssn#1 {\advance \subssscacnumber by
1{\par\medskip\noindent{\bfnn \the\sectnumber.\the\subscnumber.\the\subsscanumber.\the\subssscacnumber\ #1}\quad}}%

\newcount\subssscaenumber%
\def\caesubsssn#1 {\advance \subssscaenumber by
1{\par\medskip\noindent{\bfnn \the\sectnumber.\the\subscnumber.\the\subsscanumber.\the\subssscaenumber\ #1}\quad}}%

\newcount\subssscbanumber%
\def\cbasubsssn#1 {\advance \subssscbanumber by
1{\par\medskip\noindent{\bfnn \the\sectnumber.\the\subscnumber.\the\subsscbnumber.\the\subssscbanumber\ #1}\quad}}%

\newcount\subssscbbnumber%
\def\cbbsubsssn#1 {\advance \subssscbbnumber by
1{\par\medskip\noindent{\bfnn \the\sectnumber.\the\subscnumber.\the\subsscbnumber.\the\subssscbbnumber\ #1}\quad}}%

\newcount\subssscbcnumber%
\def\cbcsubsssn#1 {\advance \subssscbcnumber by
1{\par\medskip\noindent{\bfnn \the\sectnumber.\the\subscnumber.\the\subsscbnumber.\the\subssscbcnumber\ #1}\quad}}%

\newcount\subsssccanumber%
\def\ccasubsssn#1 {\advance \subsssccanumber by
1{\par\medskip\noindent{\bfnn \the\sectnumber.\the\subscnumber.\the\subssccnumber.\the\subsssccanumber\ #1}\quad}}%

\newcount\subsssdaanumber%
\def\daasubsssn#1 {\advance \subsssdaanumber by
1{\par\medskip\noindent{\bfnn \the\sectnumber.\the\subsdnumber.\the\subssdanumber.\the\subsssdaanumber\ #1}\quad}}%

\newcount\subsssdabnumber%
\def\dabsubsssn#1 {\advance \subsssdabnumber by
1{\par\medskip\noindent{\bfnn \the\sectnumber.\the\subsdnumber.\the\subssdanumber.\the\subsssdabnumber\ #1}\quad}}%

\newcount\subsssdacnumber%
\def\dacsubsssn#1 {\advance \subsssdacnumber by
1{\par\medskip\noindent{\bfnn \the\sectnumber.\the\subsdnumber.\the\subssdanumber.\the\subsssdacnumber\ #1}\quad}}%

\newcount\subsssdaenumber%
\def\daesubsssn#1 {\advance \subsssdaenumber by
1{\par\medskip\noindent{\bfnn \the\sectnumber.\the\subsdnumber.\the\subssdanumber.\the\subsssdaenumber\ #1}\quad}}%

\newcount\subsssdbanumber%
\def\dbasubsssn#1 {\advance \subsssdbanumber by
1{\par\medskip\noindent{\bfnn \the\sectnumber.\the\subsdnumber.\the\subssdbnumber.\the\subsssdbanumber\ #1}\quad}}%

\newcount\subsssdbbnumber%
\def\dbbsubsssn#1 {\advance \subsssdbbnumber by
1{\par\medskip\noindent{\bfnn \the\sectnumber.\the\subsdnumber.\the\subssdbnumber.\the\subsssdbbnumber\ #1}\quad}}%

\newcount\subsssdbcnumber%
\def\dbcsubsssn#1 {\advance \subsssdbcnumber by
1{\par\medskip\noindent{\bfnn \the\sectnumber.\the\subsdnumber.\the\subssdbnumber.\the\subsssdbcnumber\ #1}\quad}}%

\newcount\subsssdcanumber%
\def\dcasubsssn#1 {\advance \subsssdcanumber by
1{\par\medskip\noindent{\bfnn \the\sectnumber.\the\subsdnumber.\the\subssdcnumber.\the\subsssdcanumber\ #1}\quad}}%

\newcount\subsssdcdnumber%
\def\dcdsubsssn#1 {\advance \subsssdcdnumber by
1{\par\medskip\noindent{\bfnn \the\sectnumber.\the\subsdnumber.\the\subssdcnumber.\the\subsssdcdnumber\ #1}\quad}}%

\newcount\subsssehanumber%
\def\ehasubsssn#1 {\advance \subsssehanumber by
1{\par\medskip\noindent{\bfnn
\the\sectnumber.\the\subsenumber.\the\subssehnumber.\the\subsssehanumber\ #1}\quad}}%

\newcount\subsssfhcnumber%
\def\fhcsubsssn#1 {\advance \subsssfhcnumber by
1{\par\medskip\noindent{\bfnn
\the\sectnumber.\the\subsfnumber.\the\subssfhnumber.\the\subsssfhcnumber\ #1}\quad}}%


\chap{}{Equivalent  Inequalities}{}

\abstract{
 Equivalencies of many basic elementary inequalities  are given}

\section{Introduction} It\footnote{}{\sevenrm{\sevenit Mathematics subject classification (2000)}: 26D15.}\footnote{}{\sevenrm{\sevenit Keywords}: Equivalent inequalities geometric-arithmeti mean inequality, Bernoulli inequality, H\"older \hfill\break\vskip-.7truecm inequality, Minkowski inequality, monotonic functions, convex functions.} is an amazing act that almost all   elementary inequalities can be derived  from a few  extremely elementary 
results and are in fact equivalent to these results.  Much of this is known and is mentioned in passing in most basic books on 
inequalities, see for instance \refb{{\bf3}  pp 212--213}.  Other equivalencies arise  from the
 ability for equivalent inequalities to hide beneath almost impenetrable disguises.    In this note we collect all these various equivalencies and disguised equivalencies  so as to make them rapidly available to anyone who is interested.  The inequalities discussed here are for the most part those found in the references \refb{\bf3,8,9}. Of course there is a vast field  beyond this, see for instance \refb{\bf1,6}, and many others beyond the competence of the authors to discuss.  In this way the present paper can be regarded as a permanent work in progress as others with this competence add to the list of equivalent inequalities.

The question of equivalence can cause other problems.  While we are all agree that  inequality $\cal I$ is equivalent to inequality$\cal J$ if and only if a proof of  of $\cal J$ can be found under the assumption that $ \cal I$ holds and conversely  is not always that simple. Some equivalent inequalities are really almost equivalent  typographically: that is ${\cal I}: A\ge B $ is equivalent to ${\cal J}: B\le A$. Some arise from  not very subtle changes of notation ${\cal I}: \phi(x)\ge 0, x\ge 0 $, ${\cal J}: \phi(x^2)\ge 0,  x\in \BR$ while in others the change of variables is so intricate as to make the one an impenetrable disguise of the other and only the most perverse of  mathematicians would quote the disguise as the original.  In addition  a simple inequality may imply a more general one of which it is special case and so should perhaps not be considered as equivalent to the general form.  Finally there is the problem that many inequalities are thought of without specifying exactly all of the logical operators that are 
implicit in their  formulations.

We do not consider integral inequalities although it might be noted that where they exist they are usually seen very easily to be
 equivalent to their discrete analogues; see for instance \refb{{\bf3} pp.368--384}.  However the cases of equality and the extra 
complications are avoided here by restricting the discussion to the discrete case.

\gwide{Notations:
$\BN= \{0,1,2,\ldots\};\; \BN^*=\{1,2,\ldots\};\BN^{**}=\BN^*\setminus\{1\}=\{2,3, \ldots\};\BR$ is the set of all real numbers; $ \BR_+= \{x; x\ge0\};  \BP= \{x; x>0\}$;  $ \BP^n=( \BR_+^*)^n,\, n\in \BN^*$.\hfill\break
If $a,b\in \overline \BR,\, a\le b$  then the open and closed intervals with these endpoints are $]a,b[, [a,b]$ respectively.\hfill\break
If 
$n\in \BN^*$ then $\bvec a=(a_1,\ldots,a_n), \bvec w= (w_1,\ldots w_n)$ are $n$-tuples of positive numbers, 
that is $\bvec a, \bvec w\in  \BP^n$;  if $ a_1= \cdots = a_n$ we say the the $n$ tuple $\bvec a$ is constant; 
if $1\le k\le n$ then   $W_k= \sum_{i=1}^k w_i$. If \func{f}{\BR_+^*}{\BR_+^*} then $f(\bvec a)$ denotes the $n$-tuple $\bigl(f(a_1),\ldots.
f(a_n)\bigr)$; and $\bvec a\bvec b$ will denote the $n$-tuple $(a_1b_1, \ldots, a_nb_n)$.
$$
\displaylines
{
 \eufm G_n(\bvec a; \bvec w) = \bigl(\prod a_i^{w_i}\bigr)^{1/W_n};\cr
\eufm M_n^{[r]}(\bvec a; \bvec w) = 
\Big(\dfrac{1}{W_n}\sum_{i=1}^nw_ia_i^r\Big)^{1/r}, r\in \Bbb R^*,\cr
=\eufm G_n(\underline a; \underline w),r= 0.\cr
=\max \bvec a,\quad r= \infty\cr
=\min\bvec a,\; r= -\infty.\cr
\eufm M_n^{[1]}( \bvec a; \bvec w )= \eufm A_n(\bvec a; \bvec w )\quad \eufm M_n^{[-1]}( \bvec a; \bvec w )= \eufm H_n( \bvec a; \bvec w )\cr\eufm M_n^{[2]}( \bvec a; \bvec w )= \eufm Q_n( \bvec a; \bvec w )\cr
\eufm M_n^{[r]}( \bvec a; \bvec w )= \eufm M_n^{[r]}( a_1,\ldots a_n; w_1, \ldots , w_n )\quad etc.\cr
}
$$
If $\bvec w$ is a constant it is omitted from these notations; $ \eufm M_n^{[r]}( \bvec a)$ etc.\hfill\break
If $m\in \BN^*,\, m<n$  and $ \bvec a, \bvec w\in \BP^n$ then $\eufm A_m( \bvec a; 
\bvec w)$, etc., is taken to mean $\eufm A_m( a_1,\ldots, a_m;   w_1, \ldots, w_m)$, etc.\hfill\break
In the statement of an inequality $V$ denotes its set of validity and   $E$ denotes the set of equality; clearly $E\subset V$ and  the  equality is strict on $V\setminus E$. So formally an inequality is a triple $(V, {\cal I}, E)$\hfill\break
If $\cal I$ and $\cal J$ are two inequalities then $\cal I\equiv \cal J$ says that they are equivalent.\hfill\break
$\sim\cal I$ is the inequality 
$\cal I$ with inequality signs reversed, the reverse inequality.}
\section{Equivalent Inequalities} 
\bsubsn{Basic Notations} Given $t\in \BN^{**}, A\subseteq \BR^t$  and  \func{F}{A}{\BR} then $A=\hfill\break A_+\times A_0\times A _-$ where
$$
\displaylines
{
A_+= \{\bvec x; \bvec x\in A\and\! F(\bvec x)\ge 0\}, \quad A_0= \{\bvec x; \bvec x\in A\,\and\! F(\bvec x)= 0\},\cr  A _-= \{\bvec x; \bvec x\in A\,\and\! F(\bvec x)\le 0\},\cr
}
$$
 An {\sl inequality} $\cal I$ is a triple $\{V, E, F(\bvec x)\ge 0\}$ where: (a) $V$ is the {\sl set of validity of }\ $\cal I$, $V\subseteq A_+$; (b) $E$ is is the {\sl set of equality of } $ \cal I$,  $E\subseteq A_0\cap V$; (c) the formula $F(x)\ge 0$, often just called, by an abuse of language, the inequality. Alternatively we could write $\{V, E, F(\bvec x)\le 0\},\, V\subseteq A_-$, but this form is the same as the standard one  if $F$ is replaced by $-F$and we will not elaborate on this trivial point. Another variant is when $F= F_1-F_2$ and $\cal I$ is written $\{V, E, F_1(\bvec x)\ge F_2(\bvec x)\}$, or alternatively $\{V, E, F_1(\bvec x)\le F_2(\bvec x)\}$.

It is important to note that in general $V\subset A_+$.   The $t$ variables occurring in $F$ are of two kinds;  the basic variables, $t_1$ in number say and the parameters $t_2= t-t_1$ in number. Then we have the notation:  $\BR^t=  \BR^{t_1}\times \BR^{t_2},\; A= A_1\times A_2$  and $A_1\subseteq  \BR^{t_1}, A_2\subseteq  \BR^{t_2}$ etc. It usual to require the parameters to be such that the formula holds for all values of the variables.

\basubssn{Example} Consider $F(x,y,z, s,t) = (1-s-t)x+ sy + tz: \BP^3\times \BR^2
\mapsto \BR$; here $t_1= 3$ and $t_2= 2$. The basic geometric-arithmetic  mean inequality has $V= \BP^3\times 
T$ where $T$ is the triangle $\{(s,t); s\ge 0,\, t\ge 0,\, s+t\le 1\}$. However $F$ is positive on $T$ except at the corners so  if, as above,  $A_+= B\times C$ then $T\subset C$ but whereas the inequality holds on $T$ for all values of the variables the set $C$ depends on the variables.

In  many cases the exact value of either or both of $V$ and $E$ may be unknown 
\basubssn{Example} If $I\subset \BR$ is an interval then to  say that $ f:I\mapsto \BR$ is convex  is equivalent to the inequality $\cal I$ defined as follows. Let $ F(x,y, \lambda) = (1- \lambda)f(x) + \lambda f(y)-f\bigl((1- \lambda)x + \lambda y\bigr):I^2\times [0,1]\mapsto \BR$ then ${\cal I} = \{V, E , F(x,y, \lambda)\ge 0\}$ where $V= I^2\times [0,1]$ but in general $E$ is not known except for the obvious  $E_1=I^2\times \{0\}\cup\{1\}\subseteq E$. If $ E= \{(x,y); (x,y)\in I^2\and x=y\}\cup E_1$ then $f$ is said to be  strictly convex.

\bsubsn{Complementary and Complete inequalities} With each inequality $
\cal I$ is associated a {\sl complementary inequality}, $\sim\cal I$ being the triple $\{\sim\!V, \sim\! E, F(\bvec x)\le 0\}$  where $\sim\!V\subseteq A_-$ and $\sim\!E\subseteq A_0\cap\sim\!V$.

The pair of inequalities $\cal (I,\sim I)$ is called {\sl a complete inequality}, or more precisely {\sl the complete $\cal I$-inequality, $\tilde{\cal I}$}.
\bbsubssn{Example} If \func{F(x,y,t)}{\BP^2\times \BR}{\BR} then 
the geometric-arithmetic mean inequality, \GA,   has: $V= \BP^2\times [0,1]$ , $E= \{(x,y);x=y\}\cup\{t=0\}\cup\{t=1\}$ and the formula $F(x,y,t)\ge 0$. The complementary inequality, $\sim\!\GA$, has formula $F(x,y,t)\le 0$, $\sim\!V = \BP^2\times]\!-\!\infty, 0]\cup[1, \infty[,\, \sim\!E= E$. The complete geometric- arithmetic mean inequality, $\widetilde\GA$, is then the pair $\bigl(\GA,\sim\!\GA\bigr)$.
 \bsubsn{Equivalent Inequalities} Given two inequalities $\cal I, J$ we ay that they are equivalent, $\cal I\! \equiv\! J$, if $\cal I$ can be deduced from $\cal J$ without the use of any other inequality and vice versa. Two complete inequalities $\tilde{\cal I}, \tilde{ \cal J}$ are equivalent, $\tilde{\cal I}\! \equiv\! \tilde {\cal J}$, if both $ \cal I $ is equivalent to $ {\cal J }$ and $ \sim\!\cal I$ is equivalent to $\sim\!{\cal J}$.

There is  a minor difficulty with  this definition in that certain inequalities  are so basic that they must be allowed in any reasonable mathematical argument. These are the {\sl primitive inequalities}:

(a) inequalities between real numbers;

(b) the signs of the power functions, $x^n, n\in \BZ, x\in \BR$ or $\BR^*$. 
\bcsubssn{Examples} In all the following cases $\cal I \equivi  J$,

 (a) ${\cal I} = \{V, E, F\ge 0\};\;{\cal J}\! = \{V, E, -F\le 0\}$.

(b) ${\cal I}\! =\! \{V, E, F_1\ge F_2\}\; \;{\cal J}\! =\! \{V, E, aF_1+b\ge aF_2 +b\},\; a\in \BR_+^*,\; b\in \BR$.

(c) ${\cal I}\! =\! \{V, E, F_1\ge F_2\};\;{\cal J}\! =\! \{V, E, aF_1+b\le aF_2 +b\}.\; -a\in \BR_+^*,\;  b\in \BR$.

Clearly  this concept defines an equivalence relation in that if $
\cal I, J, K$ are inequalities then: (a) ${\cal I\equiv I}$, (b) ${\cal I \equiv  J} \Longleftrightarrow {\cal J\equiv I}$, (c) ${\cal I\equiv J \wedge J\equiv K\Longrightarrow I\equiv K}$.

As  result if we have proved ${\cal I\equiv J}$ and  ${\cal J\equiv K}$ we will not always then state that ${\cal I\equiv K}$. 
\section{The Geometric-Arithmetic Mean Inequality}
\csubsn{The Equal Weight Case}
 Perhaps the simplest  inequality is the classical one that goes back at least to the time of Euclid\noter{Euclid (fl;c.300 BC).}. 
$$
\displaylines
{
\hfill V = \{(x,y): (x,y)\in\BP^2\},\hskip 1.5cm \cr
({\cal GA}_{2,e})\hskip 3cm\sqrt{xy}\le\frac{x+y}{2}\quad \hbox{ or }\quad\eufm G_2(x,y)\le \eufm A_2(x,y): \hfill (1)\cr
\hfill    E=\bigl\{(x,y):  x=y\}.\hskip 1.5cm\cr
}
$$
 However this simple inequality is equivalent to a much more general inequality, the  equal weight geometric-arithmetic mean inequality of order $n$, where  $n\in\BN^{**}$\noter{We could of course take $ \scriptstyle n\in \BN^*$ but then $\scriptstyle n=1$ must be added to $ \scriptstyle E$; this simple  observation is often\hfill\break\vskip-.75cm valid but  will not be repeated }.
 $$
 \displaylines
 {
 \hfill V=\{(n,\bvec a); n\in\BN^{**}\land \bvec a\in  \BP^n\},\hskip 1.5cm\cr
 ({\cal GA}_{n,e})\hskip 3cm  \eufm G_n( \bvec a)\le  \eufm A_n( \bvec a):\hfill(2)\cr
\hfill E=\{\bvec a;\bvec a\in \BP^n \hbox{and $\bvec a$ is constant}\}.\hskip 1.9cm\cr
}
$$
 \casubssn{Theorem} {\sl  If $ n,m\in \BN^{**}$ then
 $$
 ({\cal GA}_{n,e}) \equiv ({\cal GA}_{m,e}).
 $$}\par
 \qed  This is immediate from known results: 
 
 (i) $\forall\; n\in \BN^*:\quad ({\cal GA}_{2,e})\Longrightarrow ({\cal GA}_{2^n,e})$; see \refb{{\bf3} pp.85-86};
 
 (ii) $\forall\; n, n'\in \BN,\, n,n'\ge 2,\;  n\ge n':\quad ({\cal GA}_{n,e})\Longrightarrow ({\cal GA}_{n',e})$; see \refb{{\bf3} p.81}.\QED
 \caasubsssn{Remark} The results used in the above theorem are due to  Cauchy\noter{Augustin Louis Cauchy,(1789-1857), a French mathematician who worked in Paris; of all mathemati-
 \hfill \vskip-.2cm cians he is the one  most often  mentioned.},  and were published in 1821.
 
 There are several different looking inequalities that are equivalent to the equal weight geometric-arithmetic mean
 inequality of order $n$; see \refb{{\bf3} pp.82--84}.
 $$
 \displaylines
 {
(a) \hfill V=\{(n, \bvec a); n\in \BN^{**} \wedge \bvec a\in \BP^n\wedge\prod_{i=1}^na_i = 1 \}, \hskip1.5cm\cr
 ({\cal I}_n)\hskip2cm n\le \sum_{i=1}^n a_i:\hfill\cr
(b) \hfill  V=\{(n, \bvec a); n\in \BN^{**} \wedge \bvec a\in \BP^n\wedge\sum_{i=1}^na_i = 1 \}, \hskip1.5cm\cr
 ({\cal J}_n)\hskip 2cm   \prod_{i=1}^n a_i \le(1/n)^n;\hfill\cr
  \hfill E=\{\bvec a; \bvec a\in \BP^n \wedge   \bvec a\; \hbox{is constant}\}.\hskip1.5cm\cr}
 $$
  \casubssn{Theorem} {\sl\ If $ n\in \BN^{**}$, then
  $$
  ({\cal I}_n)\equiv ({\cal J}_n)\equiv  ({\cal GA}_{n,e}).
  $$
  } \par
  \qed   The implications   $({\cal GA}_{n,e})\Lra({\cal I}_n)$ and $({\cal GA}_{n,e})\Lra ({\cal J}_n)$ are immediate.
  
  If now $\bvec a\in \BP^n$ and $ P=\prod_{i=1}^na_i$ the implication  $({\cal I}_n)\Lra({\cal GA}_{n,e})$ follows by applying $({\cal I}_n)$ to the $n$-tuple $\bigl((a_1/P^{1/n}),\ldots., (a_n/P^{1/n})\bigr)$.
  
  A similar argument gives the remaining implication  $({\cal J}_n)\Lra({\cal GA}_{n,e})$\QED
   \casubssn{Theorem} {\sl\  (a) Given that the logarithmic function is continuous and strictly increasing we have:
   $$
   ({\cal GA}_{2,e})\equiv \hbox{ log is strictly concave}.
   $$
   (a) Given that the exponential  function is continuous and strictly increasing we have:
   $$
   ({\cal GA}_{2,e})\equiv \hbox{exp is strictly convex}.
   $$}\par
   \qed This follows from a very simple proof in \refb{{\bf3} pp.77, 92}.\QED
  \csubsn{The General Case} After (${\cal GA}_{2,e}$) the next simplest inequality involving the geometric and arithmetic means is:
  $$
\displaylines
{
 V=V_1\times V_2,\quad V_1=\{(x,y): x,y\in \BR_+^*\},\quad V_2=\{(u,v): u,v\in \BR_+^*\},\hskip 1.5cm\cr
 ({\cal GA}_2)\hskip 1cm\;\bigl(x^uy^v\bigr)^{1/u+v}\le \frac{ux+ vy}{ u+v}\quad\hbox{or}\quad \eufm G_2(x,y;  u,v)\le\eufm A_2(x,y; u,v):\hfill(3)\cr
\hfil E=\{(x,y); x=y\}.\hskip 1.9cm \cr
 }
 $$
  Some obviously equivalent forms of this are given in the following lemma.
   \cbsubssn{Lemma} {\sl $({\cal GA}_2$) is equivalent to either of the following statements:
   $$
   \displaylines 
   {
(a)\hfill V= V_1\times V_2,\quad V_1= \{(x,y,); x,y\in \BR_+^*\},\quad V_2= \{\alpha; \alpha\in \BR \wedge 0< \alpha<1\},\cr
({\cal GA}_2)\hskip 1.5cm x^{1- \alpha}y^{ \alpha}\le (1- \alpha)x+\alpha y:\hfill(4)\cr
\hfil  \; E=\bigl\{(x,y):  x=y\bigr\}. \hskip 2.5cm.\cr
(b)\hfill V=\{(x,y,p,q); (x,y,p,q)\in \BP^4\wedge \recip{p} + \recip{q} = 1\},\hskip 1cm\cr 
({\cal Y})\hskip 2.5cm xy\le \frac{x^p}{p}+\frac{y^q}{q}:\hfill(5)\cr
\hfil    E=\bigl\{(x,y):  x^p = y^q\bigr\}.\hskip 2.5cm\cr
}
$$
}\par
\qed The first statement is an obvious rewriting of $({\cal GA}_2)$ as the inequality is just  $\eufm G_2(x,y; 1-\alpha,  \alpha)\le \eufm A_2(x,y; 1- \alpha, \alpha)$

The second is a rewriting of the first putting $ 1- \alpha = 1/p,.\, \alpha= 1/q$  and then replacing $x$ by $x^p$ and $y$ by $y^q$.\QED

Inequality {(\cal Y}) is sometimes called Young's inequality\noter{William Henry Young (1863 --1942) was one half of perhaps the most famous mathematical couple,\hfill\break\vskip -.7cm his wife being Grace Chisholm Young (1868--1944 ); they had a son who was also a famous mathemati-\hfill\break\vskip -.7cmcian Laurence Chisholm Young (1905--2000). } although it is really a very special case of that result, see \ref{{\bf5} pp.48--49}. 
\cbsubssn{Theorem} {\sl\  
$$
 ({\cal GA}_2)\equiv \forall n\in  \BN^{**} ,\; ({\cal GA}_{n,e}).
$$}\par
\qed  The one equivalence is trivial since (${\cal GA}_{2, e}$)\ is a special case of   (${\cal GA}_2$) and implies (${\cal GA}_{n,e}), \forall n\in\BN^*$ by  2.1.1 Theorem.

 The other equivalence needs more work but follows from known results; see \refb{{\bf3} pp.80--81}.
 
 (a) if $w_1. w_2$ are rational then  for a suitable $m\in \BN^*$ and $m$-tuple $\bvec b$  we can write $ \eufm A_2(\bvec a; \bvec w)$ as $\eufm A_m(\bvec b)$, and similarly for the geometric mean.  So that given  $({\cal GA}_{n,e}) \forall n\in \BN^*$  we can deduce  $({\cal GA}_2)$ when the weights are rational and get the right set $E$ since $\bvec b$ is constant exactly when $\bvec a$ is constant.
 
 (b) if  $w_1. w_2$ are real the result follows by taking the limit of the rational case except possibly for the set $E$.
 
 (c) Finally if $\bvec a$ is not constant write $w_i= q_i + r_i, i=1,2$ ,where $ q_i, i=1,2,$ is a non-zero rational.  Now
 $$
 \eqalign
 {
 \eufm A_2(\bvec a; \bvec w)=&\frac{Q_2}{W_2}\eufm A_2(\bvec a; \bvec q)+\frac{R_2}{W_2}\eufm A_2(\bvec a; \bvec r)\cr
 >&\frac{Q_2}{W_2}\eufm G_2(\bvec a; \bvec q)+\frac{R_2}{W_2}\eufm A_2(\bvec a; \bvec r),\; \hbox{by (a)}.\cr
  \ge&\frac{Q_2}{W_2}\eufm G_2(\bvec a; \bvec q)+\frac{R_2}{W_2}\eufm G_2(\bvec a; \bvec r)\; \hbox{by (b)}.\cr
   \ge&\bigl(\eufm G_2(\bvec a; \bvec q)\bigr)^{\frac{Q_2}{W_2}}\bigl(\eufm G_2(\bvec a; \bvec r)\bigr)^{\frac{R_2}{W_2}}, \; \hbox{by (b)}.\cr
   =&\eufm G_2(\bvec a; \bvec w).\cr
   }
     $$
     \QED
     \cbasubsssn{Remark} The all important last part of the proof seems to be due to Hardy, Littlewood and P\'olya\noter{Geoffrey Harold Hardy (1877--1947) and John Edensor Littlewood (1885--1977) were English and \hfill\break\vskip -.75cm George P\'olya (1887--1985) was born  in Hungary,} the trio of famous mathematicians of the first half of the twentieth century who almost single handedly founded the 
theory of inequalities with their book \refb{\bf6}. 
     
     We now generalize  3.1.1 Theorem  and agree to write: 
     $$
     \displaylines
     {
   \hfill V=\{  (n, \bvec a, \bvec w); n\in \BN^{**}, \bvec a, \bvec w\in \BP^n\},\hskip 1.5cm\cr
({\cal GA}_n) \hskip 2cm     \eufm G_n( \underline a; \underline w)\le \eufm A_n( \underline a; \underline w):\hfill(6) \cr
\hfill  E=\{\bvec a; \bvec a\in \BP^n \hbox{and $\bvec a$ is a constant}\},\hskip 1.5cm\cr    
}
 $$
     
      the geometric-arithmetic mean inequality of order $n$ .
\cbsubssn{Theorem} {\sl\  If $ n,m\in \BN^{**}$ then
$$
({\cal GA}_n)\equiv ({\cal GA}_m).
$$}\par
\qed The proof is essentially the same as that of 3.2.1 Theorem  but considerably easier.

(i) $\forall \, n\in \BN^{**}:   ({\cal GA}_n)\Lra ({\cal GA}_{n+1})$; see \refb{{\bf3} pp.90--91}.

(i) $\forall \,n\in \BN^{**}:   ({\cal GA}_{n+1})\Lra ({\cal GA}_n)$; see \refb{{\bf3}, p.81}.  This is just the Cauchy reverse induction used in the earlier theorem  As the proof in the reference is garbled  it  given here in full and in a slighty more general form.  Assume that $ n,m\in \BN^{**},\, m<n$, and let $
\bvec b$ be the $n$-tuple defined as 
$$ b_i=
\cases
{ a_i, & if $ 1\le i\le m$,\cr
 \eufm A_m( \underline a; \underline w), & if $ m< i\le n$.\cr
}
$$
 Then  by  $ ({\cal GA}_n)$,  $\eufm G_n( \underline b; \underline w)\le \eufm A_n( \underline b; \underline w)=\eufm A_m( \underline a; \underline w)$, or 
 $$
 \bigl(\eufm G_m( \underline a; \underline w)\bigr)^{W_m/W_n}\bigl( \eufm A_m( \underline a; \underline w)\bigr)^{(1-w_m)/W_n}\le  \eufm A_m( \underline a; \underline w)
 ;,
 $$
  that is $\eufm G_m( \underline a; \underline w)\le \eufm A_m( \underline a; \underline w)$. The set of equality is readily checked.  That is we have deduced $ ({\cal GA}_m)$.\QED
  \cbsubssn{The Inequalities of Rado and Popoviciu}   The Rado and Popoviciu inequalities of order $n$ are\noter{ Richard Rado (1906--1989; Tiberiu Popoviciu (1906--1975).}:
  $$
  \displaylines
  {
 \hfill V=\{(n,\bvec a, \bvec w); n\in \BN^{**} \wedge\bvec a, \bvec w \in \BP^n\},\hskip 1.5cm\cr
({\cal R}_n) \hskip 1cmW_n\bigl(\eufm A_n(\bvec a;\bvec w)- \eufm G_n( \bvec a; \bvec w)\bigr)\ge W_{n-1}\bigl(\eufm A_{n-1}( \bvec a; \bvec w)- \eufm G_{n-1}( \bvec a; \bvec w)\bigr):\hfill (7)\cr
 \hfill  E=\{\bvec a;\bvec a\in \BP^n\wedge a_n = \eufm G_{n-1}( \bvec a; \bvec w)\}.\hskip 1.5cm\cr
  \hfill V=\{(n,\bvec a, \bvec w); n\in \BN^{**} \wedge\bvec a, \bvec w \in \BP^n\},\hskip 1.5cm\cr
({\cal P}_n) \hskip 3cm \frac{\eufm A_n( \bvec a; \bvec w)}{\eufm G_n( \bvec a; \bvec w)}^{1/W_n}\ge \frac{\eufm A_{n-1}( \bvec a; \bvec w)}{ \eufm G_{n-1}( \bvec a; \bvec w)}^{1/W_{n-1}}:\hfill(8)\cr
\hfill  E=\{\bvec a;\bvec a\in \BP^n\wedge a_n = \eufm A_{n-1}( \bvec a; \bvec w)\}.\hskip 1.5cm\cr
}
  $$
   The inequality $({\cal R}_n)$ was given as an exercise in \refb{{\bf5} p.61} but has been rediscovered  many times and variants have been much studied; the multiplicative analogue $({\cal P}_n)$ was given a litle later by Popoviciu.  While $({\cal GA}_n)$, $n\in \BN^{**}$, says that if $a_1\ne a_2$ the sequence  $W_n\bigl(\eufm A_n( \underline a; \underline w)- \eufm G_n( \underline a; \underline w)\bigr),\, n\in \BN^{**}$, is positive the inequalities $({\cal R}_n), n\in \BN^{**},$ together with $({\cal GA}_2)$,  say that this sequence is positive and increasing. Although this is an apparently stronger statement the inequalities are essentially equivalent.
   
   \cbcsubsssn{Theorem} {\sl If $n\in \BN^{**}$  then:
   $$
   ({\cal R}_n)\equiv ({\cal GA}_2),\quad \hbox{and}\quad ({\cal P}_n)\equiv ({\cal GA}_2).
   $$}\par
   \qed Again this is a consequence of known results see \refb{{\bf3} p.26}. For instance 
   
  $$ W_n\bigl(\eufm A_n( \underline a; \underline w)- \eufm G_n( \underline a; \underline w)\bigr)\ge W_{n-1}\bigl(\eufm A_{n-1}( \underline a; \underline w)- \eufm G_{n-1}( \underline a; \underline w)\bigr)
  $$
   can be rewritten as
   $$
   \frac{w_n}{W_n}a_n + \frac{W_{n-1}}{W_n}\eufm G_{n-1}( \underline a; \underline w)\ge  \eufm G_n( \underline a; \underline w).
   $$
   If this is valid putting $a_n =x$ and $ a_i= y,\,, 2 \le i\le n$, gives  $({\cal GA}_2)$. On the other hand using  $({\cal GA}_2)$ proves this last inequality.\QED

\section {Bernoulli's Inequality} 
 \wide{Where necessary we put $ 0^0= 1$.}  
 \dsubsn{The Basic Bernoulli-Barrow Inequality}
    After $({\cal GA}_{2,e})$  perhaps the next simplest inequality is 
     $$
   \displaylines
  {
  V=V_1\times V_2\,\quad V_1=\{x;  x\in \BR_+\},\quad V_2=\{\alpha; \alpha\in \BR \wedge  0\le  \alpha\le 1\},\hskip 1.5cm \cr
({\cal B}_1) \hskip 3cm(1+x)^{ \alpha}\le 1 + \alpha x:\hfill(1)\cr
 \hfill E=\{(x, \alpha); x= 0\vee \alpha= 0\vee \alpha=1\}.\hskip 1cm\cr
}
  $$
 Writing $F(x, \alpha) $ for the left-hand side of  (1) and $A(x, \alpha)$ for the right-hand side and  $V= [0, \infty[\times [0,1]$ then we can express  $({\cal B}_1)$ as:
$$
\fa (x, \alpha)\in V: F(x, \alpha) \le A(x, \alpha);\qquad   E= \partial V.
$$
It is obvious that $F$ has a natural domain that is larger than $V$. The extension of the above inequality to this natural domain will be considered below. 

The  result   estimates the binomial function  $F$ by by the  linear function $A$  and   goes under the name of Bernoulli's inequality. The original Bernoulli\noter{Jacob Bernoulli (1654--1705), a Swiss mathematician who worked in
Basel.  The inequality was \hfill\break\vskip -.7cm proved 20 years earlier by the British mathematician  Isaac Barrow, (1630--1677). The inequality should \hfill\break\vskip -.7cm be called the
Barrow-Bernoulli inequality.} inequality only considered the  simplest case of the binomial function, namely $ (1+x)^n, \, n\in \BN^{**}$. 

We now consider extending the inequality  $({\cal B}_1)$ to  the natural domain of the   function $F$. 

First note an easy equivalent form of $({\cal B}_1)$:
   $$
   \displaylines
  {
 \hfill\phantom{aaaaaa}V=V_1\times V_2,\quad V_1=\{x; x\in \BR_+^*\},\quad V_2=\{\alpha; \alpha\in \BR \wedge  0\le  \alpha\le 1\},\hskip .75cm \cr
({\cal B}_1)\hskip 3cm(1+x)^{ 1-\alpha}\le 1 +(1- \alpha) x:\hfill(2)\cr
 \hfill E=\{(x, \alpha); x= 0\vee \alpha= 0\vee \alpha=1\}.\hskip2cm\cr
}
  $$
Now consider the  rest of the natural domain for the variable $x$.
$$
 \displaylines
  {
 \hfill\phantom{aa} V=V_1\times V_2,\quad\ V_1=\{x; x\in \BR\, \and\, 0<x\le 1\},\quad V_2=\{\alpha; \alpha\in \BR\, \wedge\,  0\le  \alpha\le 1\},\hskip .75cm \cr
({\cal B}_2)\hskip 3cm(1+x)^{ \alpha}\le 1 + \alpha x:\hfill(1)\cr
 \hfill E=\{(x, \alpha); x= 0\vee \alpha= 0\vee \alpha=1\}.\hskip 1cm\cr
}
 $$
\dasubssn{Theorem} {\sl\ $({\cal B}_1)\equiv  ({\cal B}_2)$}.\par
\qed (i) Consider the function $ \phi(x) = \dfrac{-x}{1+x}$; it easy to see that if $-1\le x\le 0$ then  $ \phi(x)\ge 0$, and obviously  $ \phi(x) = 0$ if and only if $x=0$. Hence by (2) 
$$
\Big(1-\frac{x}{1+x}\Bigr)^{1- \alpha}\le 1- (1-\alpha)\frac{x}{1+x},\, {\rm or}\,  (1+x)^{\alpha}\le 1+ \alpha x, 
$$
and the inequality is strict unless $x=0$, $ \alpha= 0$ or $ \alpha= 1$. 

(ii) Now  note that if  $-1\le \phi(x)\le 0$ then $x\ge 0$ and so we can reverse the preceding argument .\QED
Combining the two inequalities  (${\cal B}_1$)  and  (${\cal B}_2$)  gives
$$
 \displaylines
  {
 \hfill V= V_1\times V_2\quad V_1=\{x: x\in \BR \wedge x>-1\},\quad V_2=\{\alpha; \alpha\in \BR \wedge   0\le  \alpha\le 1\},\hskip 1.5cm \cr 
({\cal B}_3)\hskip 3cm(1+x)^{ \alpha}\le 1 + \alpha x:\hfill(1)\cr
 \hfill E=\{(x, \alpha); x= 0\vee \alpha= 0\vee \alpha=1\}.\hskip 1.5cm\cr
}
 $$
\daasubsssn{Corollary } {\sl\ $({\cal B}_3)\equiv  ({\cal B}_1)$}.\par

We finally consider the natural domain of $ \alpha$    

$$
 \displaylines
  {
 \hfill V=V_1\times V_2,\quad V_1= \{x:x\in \BR_+^*\},\quad V_2=\{\alpha;\alpha\in \BR \wedge \alpha> 1\},\hskip 1.5cm \cr
({\cal B}_4) \hskip 3cm(1+x)^{ \alpha}\ge 1 + \alpha x:\hfill(\sim1)\cr
 \hfill E=\{(x, \alpha); x= 0\}.\hskip 2.8cm\cr
}
 $$
\dasubssn{Theorem} {\sl\ $({\cal B}_1)\equiv ({\cal B}_4)$.}\par
\qed  Since  $0< 1/ \alpha<1$ we have by (1) that:
$$
(1+\alpha x)^{1/\alpha}<1+ x,\
$$
 which  gives ($\sim1$). Reversing the argument and using ($\sim1$) gives (1).\QED
 $$
 \displaylines
  {
 \hfill V=V_1\times V_2,\quad V_1= \{x:x\in \BR_+^*\},\quad V_2=\{\alpha;\alpha\in \BR \wedge \alpha<0\},\hskip 1.5cm \cr
({\cal B}_5) \hskip 3cm(1+x)^{ \alpha}\ge 1 + \alpha x:\hfill(\sim1)\cr
 \hfil E=\{(x, \alpha); x= 0\}.\hskip 1cm\cr
}
 $$
\dasubssn{Theorem} {\sl\ $({\cal B}_5)\equiv ({\cal B}_4)$.}\par
\qed
Since  $1- \alpha>1$ and   by 4.1.2  Theorem: 
$$
(1+x)^{1- \alpha}>1+ (1- \alpha)x.
$$
Using the function $\phi(x)$ introduced in the proof of 3.1.2 Theorem, 
$$
\Big(1-\frac{x}{1+x}\Bigr)^{1- \alpha}>1- (1-\alpha)\frac{x}{1+x},\, {\rm or}\,  (1+x)^{\alpha} > 1+ \alpha x, 
$$
This argument can be reversed and this completes the proof.
\QED

Now using the arguments that obtained 4.1.1 Theorem we can obtain inequalities equivalent to ${\cal B}_4$ and ${\cal B}_5$ but with the range of $x$ as in ${\cal B}_2$, that is$ -1<x\le 0$; we will neither state not prove these.

Now  state the full Bernoulli-Barrow inequality:
$$
 \displaylines
  {
\hfill V=V_1\times V_2,\quad V_1= \{x:x\in \BR\wedge x>-1\},\quad V_2=\{\alpha;\alpha\in \BR \wedge 0\le\alpha\le 1\},\hskip 1.5cm \cr
\hskip 3cm  (1+x)^{ \alpha}\le 1 + \alpha x:\hfill(1)\cr
 ({\cal B})
\hfill V=V_1\times V_2,\quad V_1= \{x:x\in \BR\wedge x>-1\},\quad V_2=\{\alpha;\alpha\in \BR \wedge \alpha<0 \vee \alpha>1\},\hskip 1.5cm \cr
\hskip 3cm(1+x)^{ \alpha}\ge 1 + \alpha x:\hfill(\sim1)\cr
 \hfill E=\{(x, \alpha); x= 0\vee \alpha= 0\vee \alpha=1\}.\hskip 1.5cm\cr
}
 $$
\dasubssn{Theorem} {\sl\ $({\cal B})\equiv ({\cal B}_1)$ }\par
\qed  Immediate from the above discussion,\QED

Also as   above we can readily write an equivalent form of  ($\cal B$) 
$$
 \displaylines
  {
\hfill V=\{(x, \alpha);  x>-1\wedge 0\le \alpha\le 1\},\hskip1.5cm\cr
\hskip 3cm  (1+x)^{1- \alpha}\le 1 + (1-\alpha) x:\hfill(2)\cr
 ({\cal B}) 
\hfill V=\{(x, \alpha); x> -1 \wedge  \alpha< 0 \, \or \alpha> 1\},\hskip1.5cm\cr
\hskip 3cm(1+x)^{1- \alpha}\ge 1 + (1-\alpha) x:\hfill(\sim2)\cr
 \hfill  E=\{(x, \alpha); x= 0\vee \alpha= 0\vee \alpha=1\}.\hskip 1.5cm\cr
}
 $$
 \dacsubsssn{} The inequality (2)  can be rewritten as follows:
 $$
 (1+x)^{ \alpha}\ge \frac{1}{1-\frac{\alpha x}{1+x}}.
 $$
 There is a similar variant of $(\sim2)$ provided $1+ (1- \alpha)x>0$; that is  if $\alpha>1$ then $-1<x<1\big/(\alpha-1)$, and if  $\alpha<0$  then $x> \big/( \alpha-1)$.

\dasubssn{The Original Bernoulli Inequality}  Before leaving  this topic it is worth noting that  the historical Bernoulli inequality is 
equivalent to ({\cal B}), as was pointed out in \refb{{\bf5} pp. 40--41}.
 $$
 \displaylines
 {
 \hfill V=V_1\times V_2,\quad V_1= \{x; x\in \BR_+^*\},\quad V_2=\{n;  n \in \BN^{**}\}\hskip 1.5cm\cr
 ({\cal OB})\hskip 3cm (1+x)^n\ge 1+nx:\hfill\cr
 \hfill E= \{x=1\}.\hskip 1.5cm\cr
 }
 $$
\daesubsssn{Theorem} {\sl\ 
$
({\cal OB})\equiv ({\cal B}).
$
}\par
\qed  The one implication is trivial and the other is in the above reference.\QED
The proof of the one implication given by Hardy, Littlewood and P\'olya is a little complicated but depends on the following interesting lemma.
\daesubsssn{Lemma} {\sl\  If $y>0$ and $n\in \BN^* $ then
$$
\frac{y^{n+1}-1}{n+1}\ge \frac{y^n-1}{n},\eqno(3)
$$
 with equality if and only if 
$ y=1$}\par
\qed  Inequality (3) is equivalent to  
$$
p(x) = ny^{n+1} -(n+1) y^n + 1\ge 0,\quad y>0.
$$
 Elementary arguments, see \refb{{\bf6} pp.3--4,} show that the polynomial $p$ has two positive roots, a double root when $y=1$.  This give (3) and the case of equality.\QED 
 \daesubsssn{Corollary} {\sl\ If $y>0$ and  $p,q\in \BN^{**},\, p>q$, then
 $$
\frac{y^p-1}{p}\ge \frac{y^q-1}{q},\eqno(3)
$$
 with equality if and only if 
$ y=1$.}\par

\dsubsn{Variants of the Bernoulli Inequality}
\dbsubssn{Changes of Variable} If  $I$ is a subinterval of $\BR$ and \func{\phi}{I}{]-1, \infty[}  then  ($\cal B$)  gives 

$$
 \displaylines
  {
\hfill V=\{(x, \alpha); \  x\in I\, \and\, 0\le \alpha\le 1\},\hskip 1.5cm:\cr
\hskip 3cm \bigl(1+\phi(x)\bigr)^{ \alpha}\le 1 + \alpha\phi(x):\hfill(4)\cr
 ({\cal I})
\hfill V=\{(x, \alpha); \  x\in I\, \and \alpha< 0 \, \or \alpha> 1\},\hskip1.5cm\cr
\hskip 3cm \bigl(1+\phi(x)\bigr)^{ \alpha}\ge 1 + \alpha\phi(x):\hfill(\sim4)\cr
 \hfill E=\{(x, \alpha); \phi(x)= 0\vee \alpha= 0\vee \alpha=1\}.\hskip 1.5cm\cr
}
 $$
  Further if $\phi$ is strictly monotonic applying $\inv{ \phi}$ we can reverse the argument and so in this case  $ (\cal I) \equiv(\cal B)$.
 
   Such an argument has been used in the proofs of 4.1.1 Theorem and 4.1.2 Theorem.
   
    Applying this idea the following examples give various  inequalities that are equivalent to ({\cal B}).

 \dbasubsssn {Example} $I = ]-\infty, 1],\;  \phi(x) = -x$ ; $ \phi$ is zero at $ x=0$. So  from  ({\cal I}):
 $$
 \displaylines
  {
\hfill V=\{(x, \alpha);   x<1\, \and\, 0\le \alpha\le 1\},\hskip 1.5cm\cr
\hskip 3cm (1-x)^{ \alpha}\le 1 - \alpha x:\hfill(5)\cr
 \hfill V=\{(x, \alpha);   x<1\, \and\,  \alpha\le 0\wedge \alpha\ge1\},\hskip 1.5cm\cr
\hskip 3cm (1-x)^{ \alpha}\ge 1 - \alpha x:\hfill(\sim5)\cr
 \hfill E=\{(x, \alpha); x= 0\vee \alpha= 0\vee \alpha=1\}.\hskip 1.5cm\cr
}
 $$
\dbasubsssn {Example}  $I= ]0,\infty[,\; \phi (x) = x-1$; $ \phi$ is zero at $ x=1$. So  from ({\cal I}):
 $$
 \displaylines
  {
\hfill V=\{(x, \alpha);  x\in \BR_+^* \and\, 0\le \alpha\le 1\},\hskip1.5cm\cr
\hskip 3cm  x^{ \alpha}\le (1 -\alpha)+ \alpha x:\hfill(6)\cr
\hfill V=\{(x, \alpha);  x\in \BR_+^* \and \alpha\le 0\vee \alpha\ge 1\},\hskip1.5cm\cr
\hskip 3cm x^{ \alpha}\ge (1- \alpha) + \alpha x:\hfill(\sim6)\cr
 \hfill E=\{(x, \alpha); x= 1\vee \alpha= 1\vee \alpha=1\};\hskip 1.5cm\cr
}
 $$
\dbsubssn{Increasing the Number of Variables} 

(A)  If $S$ is a set in $\BR^n$ and  \func{\phi}{S}{]-1,\infty[}  and  $\bvec x=(x_1, \ldots, x_n)$  then  ($\cal B$) gives 
$$
 \displaylines
  {
\hfill V=\{(\bvec x, \alpha); \bvec x\in S\, \and\, 0\le \alpha\le 1\}, \hskip 1.5cm\cr
\hskip 3cm \bigl(1+\phi(\bvec x)\bigr)^{ \alpha}\le 1 + \alpha \phi(\bvec x): \hfill(7)\cr
 ({\cal J}) 
 \hfill V\{(\bvec x, \alpha); \bvec x\in S\, \and  \alpha<0\wedge \alpha>1\},\hskip 1.5cm\cr
\hskip 3cm \bigl(1+\phi(\bvec x)\bigr)^{ \alpha}\ge 1 + \alpha \phi(\bvec x):\hfill(\sim7)\cr
 \hfill E=\{((\bvec x, \alpha): \phi(\bvec x)=0 \or \alpha= 0 \or \alpha=1\}.\hskip 1cm\cr
}
 $$
If  $ \phi(\bvec x(t))= t$ for some function $\bvec x$ then (6) implies ($\cal B$) and so  in this case ${\cal J}\equiv (\cal B)$.
\dbbsubsssn{Example}
 Let $ S=\BP^2,\,  \phi(u,v) = v/u$  and apply the above idea in the situation of 4.2.1.2 Example. Then 4.2.1.2 (5)  becomes 
 $$
\Bigl(\frac{v}{u}\Bigr)^{\alpha}\le (1- \alpha) + \alpha(v/u)\Lra u^{1- \alpha} v^{ \alpha}\le (1- \alpha)u + \alpha v,
$$
 with equality  when $ \phi(u,v) = 1$, that is if $ u=v$,
 
This however is, omitting the cases $ \alpha= 0, 1$,  just 4.2.1 (4)  
showing that  this part of the Bernoulli inequality is essentially 

equivalent to (${\cal GA}_2,$).   This can be made precise if we allow general weights in the two means.  Then we get the following extension of the geometric-arithmetic mean inequality; \refb{{\bf3} pp.148--149}. For simplicity we state it only in the form of 3.2.1 Lemma.
$$
\displaylines
{
\hfill V=\{(x,y, \alpha); x,y\in \BR_+^* \and  0< \alpha<1\},\hskip 1.5cm\cr
\hskip 3cm \eufm G_2(x,y; 1- \alpha, \alpha)\le  \eufm A_2(x,y; 1- \alpha, \alpha):\hfill;\cr
(\widetilde{\cal GA}_2)\hfill V=\{(x,y, \alpha); x,y\in \BR_+^* \and  \alpha<0  \or \alpha>1\},\hskip 1.5cm\cr
\hskip 3cm \eufm G_2(x,y; 1- \alpha, \alpha)\ge\eufm A_2(x,y; 1- \alpha, \alpha):\hfill\cr
\hfill E=\{(x,y); x=y\or \alpha= 0\or \alpha=1\}.\hskip 1cm\cr
}
$$
\dbbsubsssn{Theorem} {\sl\  $({\cal B} )\equiv (\widetilde{\cal GA}_2$).}\par
\qed Imediate from the previous discussion.\QED

\dbbsubsssn{R\"uthing's inequality}\noter{Dieter R\"uthing; but the result
 appears earlier in [{\sevenit{\bfnnn6}  pp.39--42}].};

If $S =\BP^2$,  $ \phi(a,b) = a/b-1 $;  $ \phi$ is zero  when $a=b$; and $ \phi(x,1) = x$.   From ({\cal J}) we get that
  $$
 \displaylines
  {
\hfill  V=\{(a,b, \alpha);  a,b\in \BR_+^*\wedge\; 0\le \alpha\le 1\},\hskip1.5cm \cr
 \hskip 3cm\alpha a^{ \alpha-1}(a-b)\le a^{ \alpha} -b^{ \alpha}\le \alpha b^{ \alpha-1}(a-b)\hfill \cr
 ({\cal RU})\hfill V=\{(a,b, \alpha);  a,b\in \BR_+^*\wedge\;\alpha\le 0\vee \alpha\ge 1\},\hskip 1.5cm \cr
 \hskip 3cm\alpha a^{ \alpha-1}(a-b)\ge a^{ \alpha} -b^{ \alpha}\ge \alpha b^{ \alpha-1}(a-b); \hfill(7)\cr
  \hfill E=\{(a,b); a=b\or \alpha= 0\or \alpha= 1\};\hskip 1.5cm\cr
}
 $$
 
\dbbsubsssn{Theorem} {\sl\ $({\cal RU})\equiv ({\cal B})$}\par
\qed Immediate from the general discussion.\QED

A rewriting of part of  (7):

$$
\alpha a^{ \alpha} + b{ \alpha}\ge \alpha a^{ \alpha-1}b
$$
is called Jacobsthal's inequality and being equivalent to a part of ({\cal B}) is another inequality that is equivalent to ({\cal B}),
\dbbsubsssn{Remark} It might be noted that ({\cal RU}) has a slightly different form for $V$ and $E$ in that in the left-hand inequalities $b=0$ is allowed and in the right-hand inequalities $a=0$ is allowed.

(B) Alternatively  if $S$ is a set in $\BR^n$ and if $ \psi :S\mapsto\BR$  and  $\bvec a=(a_1, \ldots, a_n)$  then  ($\cal  B$) gives 
  $$
 \displaylines
  {
\hfill V= \{(x, \bvec a); x>-1 \wedge  \bvec a\in S \wedge 0\le \psi(\bvec a)\le 1 \},\hskip 1cm\cr
\hskip 2cm (1+x)^{ \psi(\bvec a)}\le 1 + \psi(\bvec a)x: \hfill\cr
 ({\cal J}) \hfill V= \{(x, \bvec a); x>-1 \wedge  \bvec a\in S \wedge \psi(\bvec a)< 0 \, \or \psi(\bvec a)> 1\},\hskip 1cm\cr
 \hskip2cm (1+x)^{ \psi(\bvec a)}\ge 1 + \psi(\bvec a)x: \hfill\cr
 \hfill E=\{(x, \bvec a); x=0\or \psi(\bvec a)= 0\or \psi(\bvec a)= 1\}.\hskip 1cm\cr
}
 $$
 
Further if $ \psi\bigl(\bvec a(\alpha)\bigr) = \alpha$, for some choice of  the function $\bvec a$ ({\cal J}) reduces to ($\cal B$) and so is equivalent to it.  

It  is of course possible to combine various of the above  as in the following example.
\dbbsubsssn{Example}  Let  $S = \BP^2$,  $\psi(p,q) = p/q,\, 0\le p/q\le 1$ and $I= \{x; x/p>-1\},\, \phi(x) = x/p$. In this case we get  from (${\cal B}$) an inequality due to Bush, \refb{{\bf8} p.365; {\bf9} p.68}, \ref{4}.
$$
 \displaylines
  {
\hfill V=V_1\cup V_2,\; \hbox{where}\hskip 2cm\cr
\hfill V_1=\{(x,p,q); x,p,q\in \BR \wedge x\ge0\wedge 0<p\le q\vee  p\le q<-x\le0\},\cr
\hfill V_2=\{(x,p,q); x,p,q\in \BR \wedge x\le0\wedge 0\le -x<p\le q\vee  p\le q<0\},\cr
({\cal BU})\hskip 1.5cm (1+x/p)^p\le (1+ x/q)^q:\hfill(8)\cr
\hfill V= \{(x,p,q); x,p,q\in \BR \wedge p<-x\le 0<q\vee p<0\le-x<q\},\hskip .5cm\cr
\hskip 1.5cm (1+x/p)^p\ge (1+ x/q)^q:\hfill(\sim8)\cr
\hfill E=\{(x,p,q); x= 0 \vee p=q\}.\hskip 3.5cm\cr
}
$$

 Clearly (8) implies  (${\cal B}$) and so is equivalent to  (${\cal B}$). 
This result is stated in \ref{{\bf2} p.8; {\bf3} pp.36--37} but the expositions leave a lot to 
be desired. The result is a little complicated anyway but says: if  the expressions in (8) have
 
meaning and if $p<q$ then (8) holds if $ p$ and $q$ have the same sign, but $(\sim8)$ holds if $ p$ and $q$ have opposite signs.\noter{A different  proof of part of (8)  can be found in [{\sevenit{\bfnnn7}  p.365}], It is easy to get this result by a\hfill\break\vskip-.75cm direct study of the function $ \scriptstyle f(r) = (1+ x/r)^r;$  $ \scriptstyle f$  is strictly increasing on the two intervals of its domain;\hfill\break\vskip-.75cm  the domain  is the  complement of the interval $\scriptstyle ](-x)^-, (-x)^+[$; the infimum of $\scriptstyle f$ in the  left interval of the-\hfill\break\vskip-.75cm domain is the same as the supremum in the right interval, being $\scriptstyle e^x$; see below 4.3.1 Theorem.}

\dbsubssn{Inductive Variants} It is often possible to disguise an inequality such as the ones we are discussing by an apparent generalization involving induction. Thus if $ x\ge -1, \;y\ge -1 ,\, \alpha\ge 0, \beta\ge 0, 0\le
 \alpha+ \beta\le 1$, and with  no loss in generality  assume that $ y\le x$,   also assume that  not both $ \alpha$ and $ \beta$ are zero, that is $ \alpha+ \beta\ne 0$:
 $$
 \eqalignno
 {
 (1+x)^\alpha(1+y)^{ \beta}=&
 \Bigl(\bigl(1+\frac{x-y}{1+y}\bigr)^{ \alpha/( \alpha+ \beta) }(1 +y)\Bigr)^{ \alpha+ \beta}\cr
 \le &\bigl(1 + \frac{\alpha}{ \alpha+ \beta}\frac{x-y}{1+y}\bigr)^{\alpha+ \beta} (1+y)^{\alpha+ \beta}, \quad \hbox{ by ($\cal B$)},\cr
  =& \bigl(1+ \frac{\alpha}{\alpha+ \beta} x + \frac{\beta}{\alpha+ \beta} y\bigr)^{ \alpha+ \beta}\cr
  \le& ( 1+ \alpha x + \beta y)\quad \hbox{ by($\cal B$)}, & (9)\cr
 }
 $$
the inequality is strict unless either (i)  $ \alpha+\beta=0$, or (ii) $ \alpha+ \beta=1,\we x=y $,  or (iii) $ 0<\alpha+ \beta<1,\, \we  x=y=0$. Since this inequality reduces to ($\cal B$) if $ x=y$ it is equivalent to  the more elementary   result. Further this inequality is the first step in an induction to an even more general looking equivalent result, due to Pe\v cari\'c,  that is equivalent to ($\cal B$).
$$
\displaylines
{
\hfill V=\{n, \bvec a, \bvec w): n\in \BN^* \wedge a_i>-1, 1\le i\le n \wedge \bvec w\in \BR_+^n , \wedge\, 0\le W_n\le 1\} \hskip .5cm\cr
({\cal B}_n) \hskip 1cm  \prod_{i=1}^n(1+a_i)^{w_i}\le 1+\sum_{i=1}^n w_i a_i:\hfill(10)\cr
\hfill E=\{ W_n= 0\vee W_n= 1\wedge \bvec a\;\hbox{is constant} \vee 0< W_n<1\wedge \bvec a= 
\bvec0\}.\cr
}
$$
The proof of this inequality, by Pe\v cari\'c, is not readily available so will  be included for completeness.

\qed If $n=1,2$ then (10) is just $({\cal B}_1)$   and (9) respectively,  so assume that $n\ge 3$ and $({\cal B}_{n-1})$.
$$
\eqalign
{
\prod_{i=1}^n(1+a_i)^{w_i}= &\Bigl(\prod_{i=1}^{n-1}(1+a_i)^{w_i/W_{n-1}}\Bigr)^{W_{n-1}}(1+a_n)^{w_n}\cr
\le \Bigl( 1&+\recip{W_{n-1}}\sum_{i=1}^{n-1} w_i a_i\Bigr)^{W_{n-1}}(1+a_n)^{w_n},\, \hbox{by the induction hypothesis},\cr
\le 1+&\sum_{i=1}^n w_i a_i,\,\hbox{by (9), noting that $\recip{W_{n-1}}\sum_{i=1}^{n-1} w_i a_i>-1$}.\phantom{aa} \cr
} 
$$
\QED
A similar argument can be used to prove that $(\sim10)$ holds if we assume that either $w_i\le 0,\, 1\le i\le n$ or $w_i\ge 1,\, 1\le i\le n$.

\dbcsubsssn{Theorem} {\sl\  $\forall n\in \BN^*\quad ({\cal B}_n) \equiv ({\cal B})$. More generally:
\hfill
\break
\phantom{aaaaaaaaaaaa}$\forall n,m\in \BN^*\qquad ({\cal B}_n) \equiv({\cal B}_m)$. }\par

\qed  Suppose than we have $({\cal B}_n)$ for a particular $n$  then by taking $ x_1=\cdots = x_n$ we get $ ({\eufm B})$ so these inequalities 
are equivalent. More generally in this way we see that  (${\cal B}_{n_1}$) holds if and only if  (${\cal B}_{n_2}$) holds.
\QED

\dsubsn{Properties of Functions} 

The inequality ($\cal B$) can be used to obtain the properties of certain functions, that in turn are equivalent to the inequality. 
We see an example of this in 4.2.2.6 (7) which can be rephrased as saying that if $ x>0$ then the function   $f(r) = (1+x/r)^r:\BR_+^*\mapsto \BR$ is 
strictly increasing; this is elaborated n the following theorem.
\dcsubssn{Theorem}  {\sl\ If $a\in \BR, \; S=\bigl]\infty, (-a)^-\bigr[\cup\bigl](-a)^+,\infty\bigr]$ then:\hfill\break
(a) \func{f(x) = (1+a/x)^x}{S}{\BR} is strictly increasing on each of the intervals of $S$; further if $x<(-a)^-$, $y> (-a)^+$ then $f(x) > f(y)$:\hfill\break
(b) \func{g(x) =(1+a/x)^{x+a}}{S}{\BR} is strictly decreasing on each of the intervals of $S$; further if $x<(-a)^-$, $y> (-a)^+$ then $g(x) < g(y)$.}\par
\qed (a) This is just ({\cal BU}), see  4.2.2.6 note (\rn 8).

(b)  Similarly  if $ x<y$ both in one of the intervals of $S$ then $-(y-a)< -(x-a)$ and both of these quantities are  also in the same interval of $S$ and so: 
$$
\eqalign
{
\bigl(1+\frac{a}{y}\bigr)^{(y+a)} =&\bigl(1-\frac{a}{y+a}\bigr)^{-(y+a)}< \bigl(1-\frac{a}{x+a}\bigr)^{-(x+a)}, \; \hbox{by (a)}\cr
=&\bigl(1+\frac{a}{x}\bigr)^{(x+a)}.\cr
}
$$

 A similar argument can be used to complete the proof of (b).
\QED
\dcasubsssn{Corollary} {\sl\ (a) The monotonicity property of either of the  functions $f_, g$ implies that of the other.\hfill\break
(b)  The monotonicity property of each of the functions $f, g$ implies ($\cal B$).\hfill\break
(c)  All of the limits $\lim_{x\to \pm\infty}f(x)$ and $\lim_{x\to \pm\infty} g(x)$ exist and are equal to $e^a$. }\par
\qed (a) This fact for the two functions $f$ and $g$ follows  from the proof of part (c) of the theorem. The rest follows  from (b).

(b) Simple changes in variable prove this.

(c) Since $0< f< g$ it follows from the theorem   that both of the limits exist. Further, since $ 0< g(x)\big/f(x) =(1+a/x)^a$, the two limits are the same.  The value of the limit is a well known property  of the exponential function.\QED

\sectn{Mean Inequalities}  
An important part of inequalities is the area of inequalities between means. As can be seen  from 
the literature, see for instance \refb{\bf{3, 8, 9}},   there are many types of means 
but in this section we will concentrate on the power means.

The fundamental inequality between these means is the following result, called the power mean inequality\noter{
 In the equal weight case  this inequality is due to Oscar Xavier Schl\"omilch,  (1823-1901), a French-\hfill\break\vskip-.7cm born 
mathematician.
He worked in Jena and Dresden. Cauchy's techniques in analysis became well\hfill\break\vskip-.7cm  known in Germany through his textbook.
The general case was given  in 1879 by Davide Besso (1845-\hfill\break\vskip-.7cm 1906), an Italian mathematician who worked in Roma
and Modena}

$$
\displaylines
{
\hfill V= \{(n, \bvec a, \bvec w, r,s): n\in \BN^{**};\, \bvec a, \bvec w\in \BP^n;\, r,s \in \BR\,\wedge r<s  \},
\cr
({\cal R;S}_n)\hskip 1cm \eufm M_n^{[r]}(\underline a; \underline w)\le \eufm M_n^{[s]}(\underline a; \underline w):\hfill(r,s)\cr
\hfill E=\{\bvec a\in \BP^n, \bvec a\, \hbox{ is constant}\}.\hskip 4.3cm\cr
}
$$
There are several proofs of this inequality and a discussion of its history in most books on inequalities; see  \refb{{\bf3} pp.202--207}.

Special cases of this result are the geometric-arithmetic mean inequality and the harmonic-geometric
 mean inequality:
$$
\eqalignno
{
\eufm G_n(\bvec a; \bvec w)\le &\eufm A_n(\bvec a; \bvec w),&(0, 1)\cr
 \eufm H_n(\bvec a; \bvec w)\le&\eufm G_n(\bvec a; \bvec w):&(-1,0)\cr
}
$$
and $({\cal R;S}_n)$ in these cases is referred to as  $({\cal GA}_n)$  and $({\cal HG}_n)$ respectively.
 
 It follows  from this result  that  both of the limits $\lim_{r\to\pm\infty}\eufm M_n^{[r]}(\underline a; \underline w)$ exist and it can be shown that: 
 $$
 \lim_{r\to\infty}\eufm M_n^{[r]}(\underline a; \underline w)=\max \bvec a,\quad \lim_{r\to\infty}\eufm M_n^{[r]}(\underline a; \underline w)= \min\bvec a.
 $$
 In addition
 $$
  \lim_{r\to0, r\ne 0}\eufm M_n^{[r]}(\underline a; \underline w)=  \eufm G_n(\bvec a; \bvec w);
 $$
 see \refb{{\bf3} pp.175--177}.

It turns out that  many of the inequalities contained in the collection  $({\cal R;S}_n)$ are equivalent and this is what we now discuss.

\esubsn{Theorem} {\sl\ Given $n, \bvec a, \bvec w,\, r,s \in V$:\hfill\break 
(a) if  $r<s<0$ then $({\cal R;S}_n)(r,s)\equiv
 ({\cal R;S}_n)(-s, -r)$;\hfill\break
 (b) if $ r\le 0\le s$ then $({\cal R;S}_n)(r, s)\equiv ({\cal GA}_n)$;\hfill\break
 (c) if $0<r< s,$ then $({\cal R;S}_n)(r, s)\equiv({\cal R;S}_n)(1, s/r)\wedge({\cal R;S}_n)(r/s, 1)$. \hfill\break
(d) if $1<s$ then $({\cal R;S}_n)(1, s)\equiv({\cal R;S}_n)(1/s, 1)$.}
\par

\qed (a)  This follows  from the second part of the  useful identity.
If $ r, s\in \Bbb R^*,\, t=s/r,\, \underline b = \underline a^r, \underline c = \underline a^t$ then
$$
\eufm M_n^{[s]}(\underline a; \underline w) = 
\big(\eufm M_n^{[t]}(\underline b; \underline w)\big)^{t/s}
=\big(\eufm M_n^{[r]}(\underline c; \underline w)\big)^{1/t}
\eqno(1)
$$
(b) This follow from the identity ontained by taking $r=s$ in the first  identity in (1), that is:
$$
\eufm M_n^{[s]}(\underline a; \underline w)  =\big( \eufm A_n( \underline b; \underline w)\big)^{1/s}.\eqno(2)
$$
(c), (d)  These follow by again using  (2).
\QED
We now turn to equivalencies amongst the simpler inequalities (${\cal GA}_n)$. 
\esubsn{Theorem} {\sl\ (a)  (${\cal B) \equiv (\cal GA}_2)$. \hfill\break
(b) If  $m, n\in \BN^{**}$  then(${\cal GA}_{m, e})\equiv ({\cal GA}_{n,e})$.\hfill\break
(c) If  $n\in \BN^{**}$  then (${\cal GA}_{n,e})\equiv ({\cal GA}_n)$ }\par
\qed (a) See 4.2.2.2 Theorem.

(b) This is implied by the famous Cauchy backward induction, see \refb{{\bf3} pp.81--82}.

(c)  See \refb{{\bf3} pp.81--82}.\QED

%
\sectn{The H\"older and Minkowski Inequalities} 

Two other basic inequalities are the H\"older\noter{Recently the name H\"older-Rogers inequality  has been suggested as 
being more in keeping with\hfill\break\vskip-.7cm  the historical record; [{\sevenit{\bfnnn7}}]. Otto Ludwig H\"older, 
(1859-1937), German mathematician who\hfill\break\vskip-.7cm 
worked in Leipzig and proved the inequality in 1889.  Leonard James Rogers, (1862-1933), an English\hfill\break\vskip-.7cm  mathematician who worked
in Leeds and Oxford  proved an equivalent inequality in 1888. Later Frigyes \hfill\break\vskip-.7cm  Riesz (1880-1956), a
 Hungarian mathematician who
and worked in  Budapest and gave a proof for both\hfill\break\vskip-.7cm  sums and for the integrals in 1910.} and 
Minkowski\noter{{Hermann Minkowski, (1864-1909), a German mathematician who worked
in K\"onigsberg, Zurich and\hfill\break\vskip-.7cm G\"ottingen. He proved the inequality in 1896.}
} inequalities: 

$$
\displaylines
{
(a) \hfill V=\{(n,p, \bvec a, \bvec b);  n\in \BN^{**}; p\in \BR \wedge p>1;\bvec a, \bvec b\in \BP^n \},\cr
({\cal H}_{n,p})\hskip 2.3cm \sum_{i=1}^n a_ib_i\le  \Bigl(\sum_{i=1}^na_i^p\Bigr)^{1/p}\Bigl( \sum_{i=1}^nb_i^{\frac{p}{p-1}} \Bigr)^{1-1/p}:\hfill(1)\cr
E=\{\bvec a, \bvec b\in \BP^n\wedge \exists \lambda, \mu\in \BR\, \hbox{such that}\; \lambda\bvec a^p+ \mu\bvec b_i^{\frac{p}{p-1}} = \bvec 0;\cr
(b)\hfill V=\{(n,p, \bvec a, \bvec b);  n\in \BN^{**}; p\in \BR \wedge p\ge1;\bvec a, \bvec b\in \BP^n \},\cr
({\cal M}_{n,p})\hskip .6cm\Bigl( \sum_{i=1}^n (a_i+b_i)^p\Bigr)^{1/p}\le  \Bigl(\sum_{i=1}^na_i^p\Bigr)^{1/p}+\Bigl( \sum_{i=1}^nb_i^p \Bigr)^{1/p}:
\hfill(2)\cr
E=\{\bvec a, \bvec b\in \BP^n\wedge p=1  \vee p>1 \wedge \exists \lambda, \mu\in \BR\, \hbox{such that}\; \lambda\bvec a+ \mu\bvec b = \bvec 0\}\cr
}
$$

\fsubsn{Theorem} {\sl\ Given $n\in\BN^{**}, p\in \BR,\, p>1$ (${\cal H}_{n.p})\equiv({\cal M}_{n,p}$).}\par
\qed  See \refb{{\bf3}  pp,189, 212--213; {\bf10} pp.119--123}, \ref{\bf7}.\QED

The inequality (${\cal H}_{n,2}$)  is usually called the Cauchy inequality\noter{Also called the Cauchy-Schwarz-Buniakovski\u\i\ inequality
and also by any one of those three names,\hfill\break\vskip-.7cm  all namings being reasonable. Cauchy gave the inequality for finite sums in 1825.
Victor\hfill\break\vskip-.7cm  Yakovlevich Bunyakovsi\u\i,{\cyrf V.Ya. Bunyakovski\u i}, also transliterated as Buniakovsky, 
(1804-1889), a Russian\hfill\break\vskip-.7cm 
mathematician who worked  in St. Petersburg and  proved the result for integrals in 1859. Hermann \hfill\break\vskip-.7cm
Amandus Schwarz (1843-1921),  German
mathematician who worked in G\"ottingen and Berlin and pro\hfill\break\vskip-.7cmved the result in 1885 for integrals.  The inequality for spaces with 
inner product
was proved by Her- \hfill\break\vskip-.7cmmann G\"unther Grassman (1809--1877)  and  later by Hermann Klaus Hugo Weyl (1809--1955) in 1918. 
} denoted by  (${\cal C}_n$)  and we have the following result.

\fsubsn{Theorem} {\sl\  If $n\in\BN^{**}, p\in \BR,\, p>1$ then  (${\cal H}_{n,p})\equiv({\cal C}_n$).}\par

\qed  See \refb{{\bf3} pp.212--213} , \ref{\bf7}.\QED

\fsubsn{Theorem} {\sl\  If $n\in\BN^{**}, p\in \BR,\, p>1$ then (${\cal H}_{n,p})\equiv({\cal R; S}_n)(1, 1/p$).}\par
\qed See \refb{{\bf3} pp. 203, 208}.\QED
The inequality (${\cal M}_{n,2}$) is usually called the triangle inequality, written (${\cal T}_n$). 

\fsubsn{Theorem} { If $n\in\BN^{**}, p\in \BR,\, p>1$ then (${\cal M}_{n,p})\equiv({\cal T}_n$).}\par

In addition theses inequalities can be proved by induction resulting in the following equivalencies.
\fsubsn{Theorem} {\sl\  Given $p\in \BR,\, p>1$:\hfill\break (a) $({\cal H}_{2,p})\equiv \forall\, n\in \BN^{**}\; ({\cal H}_{n,p})$.\hfill\break
(b) $({\cal M}_{2,p})\equiv \forall\, n\in \BN^{**}\; ({\cal M}_{n,p})$.}\par
 \qed   See \refb{{\bf3} pp.183--185, 191--192}.\QED
The value of the parameter $p$ in both of these inequalities can be extended to all of $\BR$ and the resulting inequality is a follows.
  First let us introduce the following notation: if $p\in \BR\setminus\{1\}$ then the conjugate index $p'$ is defined by
$$
(p-1)(p'-1)=1:\; \hbox{ or: if}\;  p\ne 0,\; \recip{p} + \recip{p'} = 1,\; \hbox{ and if}\; p=0, \; p'=0.
$$
Note: $p>1\Lra p'>1$,\quad  $0<p<1 \Lra p'<0$;\quad $p<0 \Lra 0<p'<1$.
$$
\displaylines
{
\hfill V=\{(n,p, \bvec a, \bvec b);  n\in \BN^{**}; p\in \BR^*\setminus\{1\};\bvec a, \bvec b\in \BP^n \},\cr 
(\widetilde{\cal H}_{n,p})\hskip 2.3cm  \Bigl(\sum_{i=1}^n a_ib_i\Bigr)^{pp'}\le  \Bigl(\sum_{i=1}^na_i^p\Bigr)^{p'}\Bigl( \sum_{i=1}^nb_i^{p'}
\Bigr)^p:\hfill(3)\cr
E=\{\bvec a, \bvec b\in \BP^n\wedge \exists \lambda, \mu\in \BR\, \hbox{such that}\; \lambda\bvec a^p+ \mu\bvec b^{p'} = \bvec 0\
}
$$
If $p>1$ then (3) is just (1) while if $ p<1$ then (3) is the same as $(\sim1)$.

\fsubsn{Theorem} {\sl\ (a) If $0<p<1$ then $({\cal H}_{n,p})\equiv ({\cal H}_{n,1/p})$.\hfill\break
(b) (a) If $p<0$ then $({\cal H}_{n,p})\equiv ({\cal H}_{n,p'})$. }\par

\qed See \refb{{\bf6} pp.24--25}.\QED

\ffsubssn{Corollary} {\sl\ $ \forall p\in\BR^*\setminus\{1\}\; (\widetilde{\cal H}_{n,p})\equiv \forall p\in\BR,\; p>1\; ({\cal H}_{n,p})$}\par

In a similar way we can extend Minkowski's inequality.
$$
\displaylines
{
\hfill V=\{(n,p, \bvec a, \bvec b);  n\in \BN^{**}; p\in \BR^*,\, p\ge 1;\bvec a, \bvec b\in \BP^n \},\cr
\hskip 1.6cm\Bigl( \sum_{i=1}^n (a_i+b_i)^p\Bigr)^{1/p}\le  \Bigl(\sum_{i=1}^na_i^p\Bigr)^{1/p}+\Bigl( \sum_{i=1}^nb_i^p \Bigr)^{1/p}:
\hfill(2)\cr
(\widetilde{\cal M}_{n,p})\hfill V=\{(n,p, \bvec a, \bvec b);  n\in \BN^{**}; p\in \BR^*,\, p\le 1;\bvec a, \bvec b\in \BP^n \},\cr
\hskip 1.6cm\Bigl( \sum_{i=1}^n (a_i+b_i)^p\Bigr)^{1/p}\ge  \Bigl(\sum_{i=1}^na_i^p\Bigr)^{1/p}+\Bigl( \sum_{i=1}^nb_i^p \Bigr)^{1/p}:
\hfill(\sim2)\cr
E=\{\bvec a, \bvec b\in \BP^n\wedge p=1 \vee p\ne 1 \wedge\exists \lambda, \mu\in \BR\, \hbox{such that}\; \lambda\bvec a+ \mu\bvec b = \bvec 0\}\cr
}
$$
\fsubsn{Theorem} {\sl\ $ \forall p\in\BR^*, (\widetilde{\cal M}_{n,p})\equiv \forall p\in\BR,\; p\ge1\; ({\cal M}_{n,p})$}\par
\qed  See \refb{{\bf6} pp.30--32}.\QED
\fsubsn{Theorem} {\sl\ $\forall r,s\in \BR,  ({\cal R, S}_n)(r,s) \equiv \forall p\in \BR^* (\widetilde {\cal H}_{n;p})
.$\ }\par
\qed A well known proof of $({\cal H}_{n;p})$ shows that $({\cal GA}_n)\Lra  ({\cal H}_{n;p})$; see \refb{{\bf3} pp.178--179}.

It is also known that $({\cal H}_{n;p}) \Lra\forall s>1,  ({\cal R, S}_n)(1,s)$; see 6.3 Theorem.

The rest of the equivalenct then follows from   5.1 Theorem.\QED
\fsubsn{Other Equivalencies} 
It is a notorious fact that many well known 
inequalities are just  a well-known inequality  in what is often an almost 
impenetrable disguise.  This disguise is brought about by the variuous 
changes of variables described in  section 4.2 above.  We now look at some of these 
``hidden'' inequalities but first consider a cople of easy equivalencies.
\fhsubssn{Weighted Inequalities} First let us note a very elementary fact that the above inequalities are equivalent to weighted forms in which (1), (2) and (3) become:
$$
\eqalignno
{
\sum_{i=1}^n w_ia_ib_i\le & \Bigl(\sum_{i=1}^nw_ia_i^p\Bigr)^{1/p}\Bigl( \sum_{i=1}^nw_ib_i^{\frac{p}{p-1}} \Bigr)^{1-1/p},&(1w)\cr
\Bigl( \sum_{i=1}^n w_i(a_i+b_i)^p\Bigr)^{1/p}\le & \Bigl(\sum_{i=1}^nw_ia_i^p\Bigr)^{1/p}+\Bigl( \sum_{i=1}^nw_ib_i^p \Bigr)^{1/p},
&(2w)\cr
\Bigl(\sum_{i=1}^n w_ia_ib_i\Bigr)^{pp'}\le&  \Bigl(\sum_{i=1}^nw_ia_i^p\Bigr)^{p'}\Bigl( \sum_{i=1}^nw_ib_i^{p'}
\Bigr)^p;&(3w)\cr
}
$$
respectively, and where of course $\
\bvec w\in \BP^n$ .

\qed (1w) is just (1) applied to the $n$-tuples $\bvec w^{1/p}\bvec a$ and $\bvec w^{(1-1/p)}\bvec b$, and (1) is just (1w) where $\bvec w$ is constant. 
A similar argument shows that (3w) and (3) are equivalent.

(2w) is just (2) applied to the $n$-tuples $\bvec w^{1/p}\bvec a$ and $\bvec w^{1/p}\bvec b$, and (2) is just (2w) where $\bvec w$ is constant.\QED
 
Of course there is an  equivalent weighted inequality $(\sim 2w)$.
\fhsubssn{Radon's Inequality} Another inequality that is equivalent to H\"older's inequality is the following
$$
\displaylines
{
\hfill V= V_1\times V_2,\qquad V_1= \{(n, \bvec a, \bvec b); n\in \BN^{**}, \bvec a, \bvec b\in \BP^n\},\cr
\hfill V_2=\{s; s\in \BR \wedge 0<s<1\};\hskip 1.1cm\cr
\hskip2cm \sum_{i=1}^n a_i^sb_i^{1-s}\le \Bigl(\sum_{i=1}^n a_i\Bigr)^s\Bigl(\sum_{i=1}^n  b_i\Bigr)^{1-s}:\hfill(4)\cr
\hfill V= V_1\times V_3,\qquad V_3=  \{s; s\in \BR \wedge s<0 \vee s>1\};\hskip .5cm\Bigr)^{1-s}:\hfill(\sim4)\cr
\hfill E=\{\bvec a, \bvec b; \bvec a\sim\bvec b\}.\hskip 2.5cm\cr
}
$$
\qed This is seen to be equivalent to $(\widetilde{\cal H}_{n;p})$ by a simple change of variable; see \refb{{\bf3} pp.181--182}.\QED
If (4) and $(\sim4)$ are written as 
$$
\eqalignno
{
\sum_{i=1}^n \frac{a_i^s}{b_i^{s-1}}\le &\frac{\Bigl(\sum_{i=1}^n a_i\Bigr)^s}{\Bigl( \sum_{i=1}^n  b_i\Bigr)^{s-1}}&(4)\cr
\sum_{i=1}^n \frac{a_i^s}{b_i^{s-1}}\ge &\frac{\Bigl(\sum_{i=1}^n a_i\Bigr)^s}{ \Bigl(\sum_{i=1}^n  b_i\Bigr)^{s-1}}&(\sim4)\cr
}
$$
the result is called Radon's\noter{Johann Radon (1887-1956), Czech-born mathematician. He worked
in Vienna and proved the\hfill\break\vskip-.75cm inequality in 1913.} inequality.  In \refb{\bf{6} p.61} it is set as an exercise and 
in \refb{{\bf3} pp.181--182} it is Theorem 3(c) and it is not noticed to be a mere  rewriting of Theorem 3(a).

\fhsubssn{Liapunov's Inequality} The inequality known as Liapunov's\noter{Aleksandr
 Mihailovi\v c Liapunov (1857--1918),{ \cyrseven A. M. Liapunov}: also transliterated
 as Liapunoff,\hfill\break\vskip-.75cm Lyapunov a Russian mathematician who worked in Kharkov and St. Petersburg. 
He gave this inequality\hfill\break\vskip-.7cm  in 1901.} inequality is the following.
$$
\displaylines
{
\hfill V=V_1\times V_2,\qquad V_1=\{(n, \bvec x, \bvec w); n\in \BN^{**};\bvec x, \bvec w\in \BP^n\},\cr 
\hfill V_2=\{(r,s,t)\}; r,s,t\in \BR \wedge t<s<r\, \vee\, r<t<s\, \vee\, s<r<t \},\cr
\hskip1cm\Bigl(\sum_{i=1}^nw_ix_i^s\Bigr)^{r-t}\le\Bigl(\sum_{i=1}^nw_ix_i^t\Bigr)^{r-s}
\Bigl(\sum_{i=1}^nw_ix_i^r\Bigr)^{s-t}:\hfill(5)\cr
({\cal L}_{n;r,s,t})\hskip 3.8cm V=V_1\times V_3,\hfill\cr\qquad 
\hfill V_3= \{(r,s,t)\}; r,s,t\in \BR \wedge t<r<s\, \vee\, s<t<r\, \vee\, r<s<t\},\cr
\hskip1cm\Bigl(\sum_{i=1}^nw_ix_i^s\Bigr)^{r-t}\ge\Bigl(\sum_{i=1}^nw_ix_i^t\Bigr)^{r-s}
\Bigl(\sum_{i=1}^nw_ix_i^r\Bigr)^{s-t}:\hfill(\sim5)\cr
\hfill E=\{\bvec x\in \BP^n \wedge \bvec x\; \hbox{is constant}\}.\hskip 2cm\cr
}
$$

\fhcsubsssn{Theorem} {\sl\ $\forall r,s,t\in \BR\; ({\cal L}_{n;r,s,t})\equiv \forall p\in \BR\; (\widetilde{\cal H}_{n,p})$}\par

$({\cal L}_{n;r,s,t})$ follows from the weighted $(\widetilde{\cal H}_{n,p})$ by an application of the change of variables:
$$
\displaylines
{
p=\frac{r-t}{r-s}, \quad \bvec a=\bvec x^{t/p} ,\quad \bvec b=\bvec x^{r/p'}.\cr
\hbox{When}:\hfill\cr
p'= \frac{r-t}{s-t}, \quad \bvec a\bvec b= \bvec x^s.\hskip 1.7cm\cr
}
$$
 Now we see that the equal weighted  $(\widetilde{\cal H}_{n,p})$  follows from $({\cal L}_{n;r,s,t})$ by the change of variables:
$$
\displaylines
{
p=\frac{r-t}{r-s},\qquad p'= \frac{r-t}{s-t},\cr
\bvec x = \bvec a^{-1/(r-s)}\bvec b^{1/(s-t)}, \quad \bvec w = \bvec a^{r/(r-s)}\bvec b^{-t/(s-t)}.\cr
\hbox{When}:\hfill\cr
\bvec w\bvec x^s = \bvec a\bvec b,\quad \bvec w\bvec x^t = \bvec a^p,\quad \bvec w\bvec x^t=\bvec b^{p'};\;
\hbox{and }\; \bvec a^p \bvec x^{r-t}=\bvec b^{p'}.\cr 
}
$$
 Since we have noted that the weighted form of  $(\widetilde{\cal H}_{n,p})$  is equivalent to the equal weighted form 
the equivalency is proved since the cases of equality are easily discussed.\QED

 \sectn{Some Remarks} This paper originated in  an attempt to rewrite a joint paper with Yuan-Chuan Li and Cheh-Chi Yeh \noter{Yuan-Chuan Li Department of Applied Mathematics,National Chung-Hsing University, Taichung\hfill\break\vskip-.2truecm 402 Taiwan,
{\sevenit ycli@amath.nchu.edu.tw}; Cheh-Chih Yeh, Department of Information Management,Lunghwa\hfill\break\vskip-.2truecm University of Science and Technology,Kueishan Taoyuan, 
 333 Taiwan;  {\sevenit ccyeh@mail.lhu.edu.tw} or \hfill\break\vskip-.2truecm{\sevenit  chehchihyeh@yahoo.com.tw}
}.  This was  not successful but the idea for the paper came from these two  mathematicians and I recognise their efforts and inspiration. 

	The historical remarks  originate from a very intersting communication received from 
Lech Maligranda\noter{Lech Maligranda, Department of Mathematics, Lule\aa University of Technology, Lule\aa, Sweden,\hfill\break\vskip-.2truecm {\sevenit lech@sm.luth.se}} whose knowledge of the history of the various inequalities discussed far exceeds my own.

Finally it should be remarked that this paper ought to be considered a work in progress.  There are  more inequalities that should be discusssed  in the field considered here  and there are  many fields with their own inequalities.  The various equivalencies between the memebers of this vast array  has yet to be fully determined, the present paper is but  a very short introduction to the whole topic.
 \sectn{Bibliography}

[1]\ Beckenbach E F  \& Bellman R {\sl Inequalities\/},  Springer-Verlag, Berlin-Heidelberg-New York.

[2] \ P.S.Bullen,{\sl\  A Dictionary of Inequalities}, Pitman Monographs and Sureys in Pure and Applied Matheamtics 97, Addison Wesley 
 Longman, Harlow, 1998.

[3] \  P.S.Bullen, {\sl\ Handbook of Means and Their Inequalities}, Kluwer Academic Publishers, Dordrecht, 2003.

[4] \ K.A.Bush,  On an application of the mean value theorem,  {\sl Amer.\ Math.\ Monthly\/}, 62 (1955),  577--578.

[5]  A. M. Fink An Essay on the History of Inequalities, {\sl  J.\ Math.\ Anal.\ Appl.\/}, 249 (2000),  118--134.

[6] \ G.H.Hardy, J.E.Littlewood  \& G.P\'olya, {\sl\ Inequalities}, Cambridge Unversity Press, Cambridge, 1934. 

[7] \  L. Maligranda,  Equivalence of the H\"older-Rogers and Minkowski inequalities,  {\sl Math.\ Ineq.\ App.\/}, 
4 (2001),  203--207.

[8] \ D.S.Mitrinovi\'c, with P.M. Vasi\'c,  {\sl\ Analytic Inequalities}, Springer-Verlag, berlin, 1970.

[9]  \ D.S.Mitrinovi\'c, J.E.Pe\v cari\'c \& A.M.Fink, {\sl Classical and New Inequalities in Analysis}, D. Reidel, Dordrecht, 1993.
                                                                
[10] H.L.Royden. {\sl Real Analysis}. Macmillan Publ. Co., New York,1988.

\bigskip

\rightline{
 P.S.Bullen}

\rightline{Department of Mathematics}

\rightline{University of British Columbia}

\rightline{Vancouver BC} 

\rightline{Canada   V6T 1Z2}

\rightline{\tt bullen@math.ubc.ca}

\end